\documentclass{tac}

 \usepackage{amssymb, amsmath}

 \input xy
  \xyoption{all}
 
 \input diagxy

 \newtheorem{theorem}{Theorem}[section]
 \newtheorem{prop}{Proposition}[section]
  \newtheorem{lemma}[theorem]{Lemma}
  
  \numberwithin{equation}{section}

\newcommand{\la}{{\langle}}
\newcommand{\ra}{{\rangle}}

\newcommand{\mbc}{{\mathbb C}}
\newcommand{\mbf}{{\mathbb  F}}

\newcommand{\mbr}{{\mathbb R}}

\newcommand{{\tlg}}{\tilde\gamma }

\newcommand{{\tlG}}{\tilde\Gamma }
\newcommand{{\mbU}}{\mathbf U }

\newcommand{\mbbc}{{\mathbf C}}

\newcommand{\mbF}{\mathbf {F} }
\newcommand{\mbL}{\mathbf {L} }
\newcommand{\mbbb}{\mathbf {B} }
\newcommand{\mbg}{\mathbf {G} }

\newcommand{\mbk}{\mathbf {K} }

\newcommand{\mbO}{\mathbf {O} }

\newcommand{\Obj}{{\rm Obj }}
\newcommand{\End}{{\rm End }}
\newcommand{\Mor}{{\rm Mor }}
\newcommand{\Hom}{{\rm Hom }}
\newcommand{\mbV}{{\mathbf V}}
\newcommand{\mbW}{{\mathbf W}}

\newcommand{\mbC}{{\mathbf C}}

\newcommand{\mbG}{{\mathbf G}}

\newcommand{\mclD}{{\mathcal D}}

\newcommand{\mbrho}{{\mathbf \rho }}
\newcommand{\Beq}{\begin{equation}}
\newcommand{\Eeq}{\end{equation}}
\newcommand{\Beqr}{\begin{eqnarray}}
\newcommand{\Eeqr}{\end{eqnarray}}
\newcommand{\Bspl}{\begin{split}}
\newcommand{\Espl}{\end{split}}
\newcommand{{\ghat}} {(G, H, \alpha, \tau)}
\newcommand{\mcD}{{\mathcal D}}

\author{Saikat Chatterjee,  Amitabha Lahiri, 
and  Ambar N. Sengupta}
\address{Institut des Hautes \'Etudes Scientifiques \\
35 Route de Chartres\\
Bures-sur-Yvette-91440
France\\ [5pt]
S.~N.~Bose National Centre for Basic Sciences \\ Block JD,
  Sector III, Salt Lake, Kolkata 700098 \\
  West Bengal, India\\[5pt]
Department of Mathematics\\
  Louisiana State University\\  Baton
Rouge, Louisiana 70803, USA
}

\thanks{  }

\title {Twisted actions of categorical groups}

\copyrightyear{2014}

\keywords{Representations; Double Categories; Categorical Groups; 2-Groups; Semidirect Products}
\amsclass{Primary18D05; Secondary: 20C99}

\eaddress{saikat.chat01@gmail.com\CR 
 amitabhalahiri@gmail.com \CR ambarnsg@gmail.com }

\newtheorem{theorem}{Theorem}

\newtheoremrm{rem}{Remark}

\let\pf\proof
\let\epf\endproof
\mathrmdef{Hom}
\mathbfdef{Set}

\begin{document}

\maketitle
\begin{abstract}
  We develop a  theory of twisted actions of  categorical groups  using a notion of semidirect product of categories. We work through numerous examples to demonstrate the power of these notions. Turning to representations, which are actions that respect vector space structures,  we establish  an analog of Schur's lemma in this context.   Keeping new terminology to a minumum, we concentrate on examples exploring the essential new notions introduced. \end{abstract}


\section{Introduction}\label{intro}

A categorical group is comprised of a pair of groups, one of which serves as the set of objects of a category and the other as the set of morphisms. Categorical groups reflect symmetries of geometric structures that also have fruitful formulations in the language of objects and morphisms. Such geometric structures are the subject of study in a large and growing body of literature  \cite{AbbasWag, Baez, BaezSchr,   BaezWise, Bart,CLSpath, CLScb,Pfeif,MarPick,SW1, SW2,Vien}, to name a few.
A natural direction of inquiry is the theory of representations of these structures; works in this direction include   \cite{BaezBFW, BarrMack, CS, Elgu, Yet}.  However,  the theory is at an early enough stage that even what should constitute a representation of a categorical group is not definitively settled.

 In this paper we  proceed in a flexible way, focusing more on {\em actions} of categorical groups and exploring numerous examples, many of which involve some vector space structure as well.  We introduce a key procedure: forming a {\em semidirect product} of  two given categories using a certain {\em twisting map}  and we present numerous examples of this construction. We then {\em explore actions of categorical groups that involve a twisting} in the way morphisms are composed. Again we work through several examples, devoting   section  \ref{S:P} entirely to a detailed development of an example that is motivated by the physical context of fields transforming under the Poincar\'e group. Finally, turning to the narrowest sense of a representation of a categorical group on a categorical vector space,  we prove results showing that such representations have very special properties that uniquely determine them from their behavior on some subspaces.

A categorical group $\mbg$, as we explain in greater detail in section \ref{S:catgrpsumm}, is comprised of an object group  $\Obj(\mbg)$, a morphism group $\Mor(\mbg)$, along with source and target homomorphisms $s, t:\Mor(\mbg)\to\Obj(\mbg)$. Thus in $\Mor(\mbg)$ there is both a {\em multiplication}  operation $(k_2,k_1)\mapsto k_2k_1$ and a {\em composition} of morphisms  $(k',k)\mapsto k'\circ k$ (when defined). In our point of view an {\em action}     of $\mbg$ on a category $\mbbc$ should involve both an action of $\Obj(\mbg)$ on $\Obj(\mbbc)$ and an action of $\Mor(\mbg)$ on $\Mor(\mbbc)$.  Thus, in particular, 
\begin{equation}
(k_2k_1)f=k_2(k_1f)
\end{equation}
for all $k_1, k_2\in \Mor(\mbg)$ and $f\in\Mor(\mbbc)$. One might  further specify a behavior of the compositions
\begin{equation}\label{E:k1f1comp}
(k_2f_2)\circ (k_1f_1),
\end{equation}
for $f_1, f_2\in\Mor(\mbbc)$. The simplest such specification would be to require that this be equal to
\begin{equation}\label{E:compsemid0}
(k_2\circ k_1)(f_2\circ f_1).
\end{equation}
A richer extension of this idea is to require that the composition (\ref{E:k1f1comp}) be equal to
\begin{equation}\label{E:compsemid}
\bigl(\eta(k_2,f_1)\circ k_1\bigr)(f_2\circ f_1),
\end{equation}
where $\eta(k_2,f_1)$ {\em reflects a twisting of $k_2$ with respect to} $f_1$. Alternatively, one could require that the twisting act on $f_2$, with  (\ref{E:k1f1comp}) being equal to
\begin{equation}\label{E:compsemid2}
(k_2\circ k_1)\bigl(\eta(k_1,f_2)\circ f_1).
\end{equation}
We will explore the choices (\ref{E:compsemid0}) and (\ref{E:compsemid}). The $\eta$-twist in (\ref{E:compsemid}) is   easier to understand when the $k$'s are 
written on the right:
\begin{equation}\label{E:compsemid3}
(f_2k_2)\circ (f_1 k_1)= (f_2\circ f_1)\bigl(\eta(k_2,f_1)\circ k_1\bigr),
\end{equation}
where the $\eta$-twist shows up as a price to pay for moving $f_1$ to the left past $k_2$, there being no `price' in the case (\ref{E:compsemid0}). For consistency of conventions we write the $k$'s to the left of the $f$'s and then (\ref{E:compsemid3}) becomes  (\ref{E:compsemid}).  This simple extension of (\ref{E:compsemid0}) leads to some remarkably rich examples. 

We begin in section \ref{S:catgrpsumm} with a summary of categorical groups, including some examples of interest for  us. In section  \ref{s:rep} we introduce the  notion of {\em action of a categorical group on a category} (following up on the requirement (\ref{E:compsemid0}) mentioned above);   no vector space structure is involved in this. We show how a double category arises from such an action. In section \ref{s:matrep} we use crossed modules, which are more concrete renditions of categorical groups, to explore an example of an action arising from a traditional representation of the object group on a vector space.
Next, in section \ref{s:semidir}, we introduce a central notion of this work: the {\em semidirect product of two categories} obtained by using a map $\eta$ (motivated by (\ref{E:compsemid}) above).  We present several examples in detail.    In section \ref{s:semirepre} we turn to the study of {\em twisted actions}, these being actions that satisfy the $\eta$-twist condition (\ref{E:compsemid}) for  the behavior of composition; we prove a double category result. Section \ref{S:P} is devoted to a detailed study of a complex example that involves more than one  twist and is motivated by the context of fields transforming under a symmetry group that is a semidirect product (such as the Poincar\'e group). In section \ref{s:irred} we study {\em twisted representations}, by which we mean twisted actions on {\em vector categories}, these being categories that involve some vector space  structure  and linearity; we  introduce a notion of {\em irreducibility }for  such representations in this context and {\em prove an analog of Schur's Lemma} (Theorem \ref{T:schur}). Finally, in section  \ref{S:exm} we study representations on {\em categorical vector spaces}, these being the vector space analogs of categorical groups, and show that they  have very special properties.

In the course of developing our theory we obtain several categories, constructed from familiar settings involving groups or function spaces, that have interesting structure and where even the fact that they form categories is sometimes unexpected. We have kept new terminology to a minimum, focusing on the essential new ideas which we develop through detailed examples.

\section{Categorical groups in summary}\label{S:catgrpsumm}

In this section we summarize some  essential facts about categorical groups. The reader may consult our earlier work  \cite[sec. 4] {CLScb} for a more systematic presentation that includes proofs. Other works include   \cite{Awody, {BaezSchr},  BarrMack,  Bart,   ForBar,{KPT}, Kelly, Macl, {Port },Whit}.

By a {\em categorical group} we   mean a category $\mbk$ along with a bifunctor 
\begin{equation}\label{E:Kotimes}\otimes:\mbk\times\mbk\to\mbk \end{equation}
such that both $\Obj(\mbk)$ and $\Mor(\mbk)$ are groups under the operaton $\otimes$ (on objects and on morphisms, respectively). The term {\em strict $2$-group} is also widely used
instead of categorical group.

Sometimes we shall write $ab$ instead of $a\otimes b$. The  functor  $\otimes$ carries the identity morphism $(1_a,1_b)$, where $1_x:x\to x$ denotes the identity morphism at $x$,  to the identity morphism $1_{ab}$:
$$1_a  1_b=1_{a  b}.$$
Taking $a$ to be the identity element $e$   in $\Obj(\mbk)$, it follows that 
\begin{equation}\label{E:1eid}
\hbox{$1_e$ is the identity element in the group $\Mor(\mbk)$. }
\end{equation}
As a very special case of this we see that if $\Obj(\mbk)$ and $\Mor(\mbk)$ are additive abelian groups then 
\begin{equation}\label{E:100}
1_0=0,
\end{equation}
where $0$ on the left is in $\Obj(\mbV)$ and the $0$ on the right is in $\Mor(\mbV)$.

The functoriality of $\otimes$ implies the `exchange law':
\begin{equation}\label{E:exchange}
(g'  f')\circ(g  f)=(g'\circ g)  (f'\circ f)\end{equation}
for all $f,f',g,g'\in\Mor(\mbk)$ for which the composites $g'\circ g$ and $f'\circ f$ are defined.

Furthermore,  for any morphisms $f:a\to b$ and $h:b\to c$ in $\mbk$, the composition $h\circ f$ can be expressed in terms of the product operation (written as juxtaposition) in $\mbk$:
\begin{equation}\label{E:fhcompprod} 
h\circ f =f1_{b^{-1}}h=h1_{b^{-1}}f. \end{equation}
For a proof see, for example,    \cite[Prop. 4.1]{CLScb}. 
In particular,  
\begin{equation}\label{E:hkkh}
hk=h\circ k=kh\qquad\hbox{if $t(k)=s(h)=e$.} \end{equation}

 There is a useful alternative formulation  \cite[Prop. 4.2]{CLScb}  of the notion of categorical group. If $\mbk$, with operation $\otimes$, is a categorical group then  the source and target maps
$$s, t: \Mor(\mbk)\to \Obj(\mbk)$$
are group homomorphisms, and so is the identity-assigning map
$$ \Obj(\mbk)\to\Mor(\mbk): x\mapsto 1_x.$$
Conversely, if $\mbk$ is a category for which both $\Obj(\mbk)$ and $\Mor(\mbk)$ are groups,  $s$, $t$, and $x\mapsto 1_x$ are homomorphisms, and the exchange law (\ref{E:exchange}) holds then $\mbk$ is a categorical group. 

An even more concrete formulation is obtained by realizing that for a categorical group $\mbg$,  the morphism group $\Mor(\mbg)$  is isomorphic to a semidirect product:
\begin{equation}\label{E:morcatcross}\Mor(\mbg)\simeq H\rtimes_{\alpha}G:\phi\mapsto \left(\phi 1_{s(\phi)^{-1}}, s(\phi)\right),
\end{equation}
where
\begin{equation}\label{E:catcross}
\hbox{$G=\Obj(\mbg)$ and $H=\ker\left(s:\Mor(\mbg)\to\Obj(\mbg)\right)$,}
\end{equation}
and
$$\alpha:G\to{\rm Aut}(H)$$
is the homomorphism given by
\begin{equation}\label{E:alphaconj}
\alpha(g)(h)=1_{g}h1_{g^{-1}}.
\end{equation}
The target map $t:\Mor(\mbg)\to\Obj(\mbg)$ is given by
\begin{equation}\label{E:thgtau}
t(h,g)=\tau(h)g,
\end{equation}
where $\tau:H\to G$ is a homomorphism, corresponding to the restriction of $t$ to $\ker s$. 
Multiplication in $\Mor(\mbg)$ corresponds to the semidirect product  group product on $H\rtimes_{\alpha}G$
given by:
\begin{equation}\label{E:crossedprod}
(h_2,g_2)   (h_1,g_1)=\bigl(h_2\alpha(g_2)(h_1), g_2g_1\bigr) \end{equation}
and  composition corresponds to the operation
\begin{equation}\label{E:crossedcomp}
(h',g')\circ (h, g) = (h'h, g), \quad\hbox{if $g'=\tau(h)g$.}
\end{equation}
The homomorphisms $\alpha$ and  
$$\tau:H\to G: \phi\mapsto t(\phi),$$
satisfy the Peiffer identities \cite{Peif}:
\begin{eqnarray}\label{E:Peiffer} 
\tau\bigl(\alpha(g)h\bigr)& =& g   \tau(h)  g^{-1}\\
\alpha\bigl(\tau(h)\bigr)(h') &=& hh'h^{-1}\nonumber
\end{eqnarray} 
for all $g\in G$ and $h\in H$.  These conditions ensure that the crucial exchange law (\ref{E:exchange}) holds. A system 
$$(G, H, \alpha,\tau),$$ 
where $G$ and $H$ are groups, $\alpha:G\to {\rm Aut}(H)$ and $\tau:H\to G$ are homomorphisms satisfying (\ref{E:Peiffer}), is called a {\em crossed module} \cite[sec 2.9]{Whit} and   \cite{Peif}.    Conversely, a crossed module $(G, H, \alpha,\tau) $ gives rise to a categorical group $\mbg$ with 
  \begin{equation}\label{E:GHcrossed}
  \hbox{$\Obj(\mbg)=G$ and $\Mor(\mbg)=H\rtimes_{\alpha}G$,}
  \end{equation}
  and
  \begin{equation}
  \hbox{ $s(h,g)=g$ and $t(h,g)=\tau(h)g$ for all $h\in H$ and $g\in G$.}
  \end{equation}
   
When working with   Lie groups, it is  natural to  require that the source, targets, and identity-assigning maps are smooth. Thus    a {\em categorical Lie group } is  a category $\mbk$ along with a functor $\otimes$ as above, such that $\Mor(\mbk)$ and $\Obj(\mbk)$ are Lie groups, and the maps $s$, $t$, and $x\mapsto 1_x$ are smooth homomorphisms. The corresponding crossed module involves Lie groups and smooth homomorphisms.

{\bf Example CG1.} For any group $K$ we can construct a categorical group $\mbk_0$ by taking $K$ as the object set and requiring there be a unique morphism $a\to b$ for every $a, b\in K$. At the other extreme we have the discrete categorical group $\mbk_d$ whose object set is $K$ but whose morphisms are just the identity morphisms. ${ }$
 
   {\bf Example CG2.} Let $H$ be an abelian group,  $G$ any group, and $\alpha:G\to {\rm Aut}(H)$ is a homomorphism. Then a categorical group is obtained using the crossed module $(G, H, \alpha,\tau)$ where $\tau:H\to G:v\mapsto e$ is the constant map carrying all elements of $H$ to the identity element of $G$. A useful if somewhat trivial case is when $G$ is the trivial group. A case of  greater interest is obtained from the group of affine automorphisms of a vector space $V$. For this we take $G$ to be $GL(V)$, the usual group of invertible linear maps $V\to V$, and $H$ to be $V$ itself, viewed as the group of translations on $V$. We take
   $$\alpha:G= GL(V)\to {\rm Aut}(V) $$
   to be the inclusion map in the sense that for every $g\in G$ we take $\alpha(g):V\to V$ to be  the mapping $g$ itself (with $V$ being viewed just as an abelian group under addition, the elements of  ${\rm Aut}(V)$  are the addition-preserving bijections $V\to V$ and so, for a general field of scalars,  an element of ${\rm Aut}(V)$ might not preserve multiplication by scalars). 
   Then the group of morphisms of this categorical group is
   $$V\rtimes_{\alpha} GL(V),$$
   which is the group of affine automorphisms of the vector space $V$. In the case where the field is $\mbr$ or $\mbc$, and $V$ is equipped with a metric, we also have the categorical group for which $H=V$ and $G$ is the orthogonal or unitary group for $V$.  More generally, let $R$ be a representation of a group $G$ on a vector space $V$. Taking $\alpha:G\to {\rm Aut}(V)$ to be given by $R$ we obtain a categorical group whose object group is $G$ and morphism group is $V\rtimes_R G$. Since $\tau$ is trivial, the source and target maps coincide:
   \begin{equation}\label{E:sthgg}
   s(h,g)=t(h,g)=g.
   \end{equation}
   We note that the group {\em product} on the morphisms is given by
   \begin{equation}\label{E:semidprodgrp}
   (v_2,g_2)(v_1,g_1)=\bigl(v_2+g_2v_1, g_2g_1)
   \end{equation}
   which is the same as composition of the affine transforms of $V$ given by $x\mapsto g_iv+v_i$ (in other words, the group product on the morphisms is the same as the group product in ${\rm Aut}(V)$). Categorical groups of this  type have been used widely in the literature; see for example  \cite{BaezBFW, CS, Pfeif}.

\section{Actions and Double Categories}\label{s:rep}

A representation of a group $G$ on   a finite-dimensional vector space $V$ is given by a mapping
$$\rho:K\times V \longrightarrow V : (k,v)\mapsto \rho(k,v)=\rho(k)v,$$
such that each $\rho(k)\in {\rm End}(V)$, and  $\rho(kk')=\rho(k)\rho(k')$  for all $k, k' \in K$.    On what type of object should a categorical group be represented? A natural choice would be a categorical vector space (defined analogously to a categorical group). We will study such representations  in section  \ref{S:exm} and show that they have very special features. 

Let us note also a more abstract view of representations.  Instead of the group $K$ let us consider the category $\mbg(K)$, which has just one object  and the morphisms are the elements of the group $K$, with composition being given by the binary group operation. Instead of the vector space $V$ let us consider a category ${\rm {\bf Vec}}(V)$ which, again, has just one object, and whose morphisms are the endomorphisms of $V$, with composition being given by the usual composition of endomorphisms. Then the representation $\rho$ of $K$ on $V$ is equally well described by a functor
$$\mbg(K)\to {\rm {\bf Vec}}(V).$$
An extension of this point of  view leads to a notion of representations of $2$-groups on $2$-categories; for this see, for example,   \cite{BaezBFW}. 
Our approach will be quite different.

 By an {\em action}  of a  categorical group $\mbG$ 
  on a  category  $\mbV$ we mean a functor
\begin{equation}\mbrho:\mbg \times \mbV \longrightarrow \mbV, \label{def:repfunc}\end{equation}
such that 
\begin{equation}\label{E:rhokkprf} 
\rho\left(k', \rho(k,f)\right)=\rho(k'k,f)
\end{equation}
for all $k,k'\in \Mor(\mbG)$ and all $f\in\Mor(\mbV)$, and
\begin{equation}\label{E:rho1ef}
\rho(1_e,f)=f
\end{equation}
for all $f\in\Mor(\mbV)$, with $1_e\in\Mor(\mbG)$ being the identity morphism of the identity element $e\in\Obj(\mbG)$.  We will often write the effect of applying $\rho$ simply by juxtaposition:
\begin{equation}\label{E:rhoax}
ax\stackrel{\rm def}=\rho(a,x).
\end{equation}

Let us now unravel the details of the structure that lies behind the compact description of an action we have just seen.

Using the fact that $\rho$ is a functor we have  
$$\rho\left(s(1_e),s(1_v)\right)=s\left(\rho(1_e,1_v)\right)=s(1_v) \qquad\hbox{for all $v\in\Obj(\mbV)$,}$$
which means
\begin{equation}\label{E:rhoevv}
\rho(e,v)=v \qquad\hbox{for all $v\in\Obj(\mbV)$.}
\end{equation}
Taking $k'=1_a$, $k=1_b$ and $f=1_v$, where $a,b\in\Obj(\mbg)$ and $v\in\Obj(\mbV)$, in (\ref{E:rhokkprf}) we have
$$\rho\left(1_a,\rho(1_b,1_v)\right)=\rho(1_{ab},1_v),$$
and then considering the sources of both sides we conclude that
\begin{equation}\label{E:rhoabv}
\rho\left(a,\rho(b,v)\right)=\rho(ab,v).
\end{equation}
Thus {\em  the conditions on $\rho$ imply that $\rho$ gives an action of the group $\Mor(\mbg)$ on the set $\Mor(\mbV)$ and also of the group $\Obj(\mbg)$ on the set $\Obj(\mbV)$.}

The conditions on $\rho$ go well beyond simply the requirement that it specify actions of  the groups $\Obj(\mbg)$ and $\Mor(\mbg)$. That $\rho$ is a functor also says something about the interaction of compositions of morphisms:
\begin{equation}\label{E:rhofunctf12k12}
\rho(k_2\circ k_1, f_2\circ f_1)=\rho(k_2,f_2)\circ\rho(k_1,f_1)
\end{equation}
for all $k_2, k_1\in\Mor(\mbg)$ for which the composition $k_2\circ k_1$ is defined and all $f_2, f_1\in\Mor(\mbV)$ for which $f_2\circ f_1$ is defined. Stated in simpler notation this reads
\begin{equation}\label{E:rhofunctf12k12b}
 (k_2\circ k_1)(f_2\circ f_1)= (k_2 f_2)\circ (k_1 f_1).
\end{equation}

A special case of interest is when $\mbV$ is a  category whose object set  is equipped with a vector space structure.  In this case we  might be interested in the case where  the map
$$\rho(a,\cdot):\Obj(\mbV)\to\Obj(\mbV): v\mapsto\rho(a,v)$$
is linear for each $a\in\Obj(\mbg)$. When $\Mor(\mbV)$ is also equipped with a vector space structure the corresponding condition would be linearity of 
$$\rho(k,\cdot):\Mor(\mbV)\to\Mor(\mbV): f\mapsto\rho(k,f)$$
 for all $k\in\Mor(\mbG)$. Needless to say other structures may also be incorporated; for example, if $\Obj(\mbg)$ is a Lie group and $\Obj(\mbV)$ is a Hilbert space we might require that $\rho$ be continuous.

\noindent{\bf Example RG1.} Let $G$ be any group,  and let  $\mbg_0$ be the categorical group with object set $G$, with a unique morphism $a\to b$, written as $(b,a)$,  for every $a, b\in G$. Let $V$ be a vector space and $\mbV_0$   the category whose object set is $V$, with a unique arrow from $v_1$ to $v_2$, for every $v_1, v_2\in V$. The law of composition is the natural one: $(v_3,v_2)\circ (v_2, v_1) =(v_3,v_1)$.  We note that $\Mor\bigl( \mbV \bigr)\simeq V^2$ is a vector space. Now if $\rho$ is a representation of $G$ on $V$, then it defines a representation of    $\mbg_0$ on   $\mbV_0$   by 
\begin{eqnarray}\label{E:RG1}
 (a, v) & \mapsto & \rho(a, v)\\
 \bigl((b,a);(v_2, v_1)\bigr) &\mapsto & \bigl(\rho(b, v_2),  \rho(a, v_1)\bigr),\nonumber
 \end{eqnarray} 
 for any $a, b\in G$ and $v_1, v_2 \in V$. ${ }$

Our objective in the remainder of this section is to  show that an action of a categorical
group gives a `double category' over  $\mbV$.  Terminology and description of this notion, originally introduced by Ehresmann, varies in the literature; a very compact, if somewhat opaque, definition is provided in    \cite[p. 44]{Macl}: a double category is a set  that is the set of morphisms of two categories such that the two composition laws obey an exchange law of the form given below in (\ref{E:exchDC}).

By a {\em double category } $\mbC_{(2)}$ {\em over a category} $\mbC$ we mean a category for which: 
\begin{itemize}
\item[(a)] the objects of $\mbC_{(2)}$ are the morphisms of $\mbC$:
$$\Obj(\mbC_{(2)})=\Mor(\mbC);$$
 \item[(b)] the morphisms of $\mbC_{(2)}$ are the  arrows of a second category $\mbC_{(2)}^h$:
$$ \Mor(\mbC_{(2)})= \Mor(\mbC_{(2)}^h);$$
  and 
  \item[(c)] the composition operation  in the category $\mbC_{(2)}^h$ and the composition operation in the category    $\mbC_{(2)}$ satisfy the conditions (h1) and (h2) stated below. 
  \end{itemize}
  Before stating conditions (h1) and (h2) we adopt some notational and terminological conventions: 
\begin{itemize}
\item[(i)]  a morphism of $\mbC_{(2)}$ will be displayed as a `vertical' arrow\begin{equation} \label{E:vertmor}
 \xymatrix{
        u_1 \ar@{-}[r] & s(F)\ar[d]_{F}\ar[r] & u_2  
        \\
   u'_1 \ar@{-}[r] & {t(F)}\ar[r] & u'_2, 
   }  
   \end{equation}
   where  the source $s(F)$ and the target $t(F)$, both being objects of $\mbC_{(2)}$, are morphisms of $\mbC$;
   \item[(ii)] the composition of morphisms in   $\mbC_{(2)}$ will be  denoted
$$(G,F)\mapsto G\circ  F  $$ 
and may be called {\em vertical  composition} in light of the display convention in (i);
\item[(iii)] the composition of morphisms in   $\mbC_{(2)}^h$ will be called {\em horizontal composition} and denoted
$$(G,F)\mapsto G\circ_h F. $$ 
\end{itemize}
 The two composition  laws, vertical and horizontal, are required to satisfy the following conditions:
\begin{itemize}
\item[(h1)] The horizontal composite $G\circ_h F$ is defined if and only if both the composites  $s(G)\circ s(F)$ and $t(G)\circ t(F)$ are defined in   $\Mor(\mbC)$ and then
the source and target of $G\circ_hF$ in $\mbC_{(2)}$ are:
\begin{equation}\label{E:stdccomp}
\hbox{ $s(G\circ_h  F)=s(G)\circ s(F)$ and $t(G\circ_h F)=t(G)\circ t(F)$,}
\end{equation}
as illustrated by the diagram:

\begin{equation} \label{E:horcomp}
 \xymatrix{
        u_1 \ar@{-}[r] & s(F)\ar[d]_{F}\ar[r] & u_2  \ar@{-}[r]& s(G)\ar[d]_{G} \ar[r]& u_3 &  \ar@{}[d]_{=}&   u_1 \ar@{-}[r] & s(G\circ_hF)\ar[d]_{G\circ_hF}\ar[r] &  u_3 
        \\
   u'_1 \ar@{-}[r] & {t(F)}\ar[r] & u'_2  \ar@{-}[r]&  t(G)  \ar[r]& u'_3  &&   
     u'_1 \ar@{-}[r] & t(G\circ_hF) \ar[r] &  u'_3 
   }  
   \end{equation} 
\item[(h2)] The exchange law
\begin{equation}\label{E:exchDC}
(G'\circ G)\circ_h (F'\circ F)=(G'\circ_h F')\circ (G\circ_h F) 
\end{equation}
holds in the diagram 
\begin{equation} \label{E:exchangelaw}
 \xymatrix{
        u_1 \ar@{-}[r] &f\ar[d]_{F}\ar[r] & u_2  \ar@{-}[r]& g\ar[d]_{G} \ar[r]& u_3         \\
   u'_1 \ar@{-}[r] & {f'}\ar[d]_{F'}\ar[r] & u'_2  \ar@{-}[r]& g'\ar[d]_{G'}  \ar[r]& u'_3\\
   u''_1 \ar@{-}[r] & {f''}\ar[r] & u''_2  \ar@{-}[r]&  g''  \ar[r]& u''_3     }  
   \end{equation} 
   for all morphisms $F$, $G$, $F'$, and $G'$  in $\mbC_{(2)}$ for which the sources and targets are related as illustrated.

\end{itemize}

\begin{prop}\label{prop:doublerep}
An action $\rho$ of a categorical group $\mbG$  on a   category $\mbV$ yields a double category over $\mbV$ whose morphisms are the elements
$$(k,f)\in \Mor(\mbG)\times \Mor(\mbV),$$
with source and target being given by:
\begin{equation}\label{E:stMorGV}
s(k,f)=f\qquad \hbox{and}\qquad t(k,f)= \rho(k)f.
\end{equation}
 Composition  of morphisms is given by
 \begin{equation}\label{E:compMorGV}
(k',f')\circ (k,f)= (k' k, f),
\end{equation}
whenever $f'=kf$, 
and horizontal composition is given by
\begin{equation}\label{E:horcompMorGV}
(k_2,f_2)\circ_h (k_1, f_1)= (k_2\circ k_1, f_2\circ f_1),
\end{equation}
whenever the compositions on the right are defined.
  \end{prop}
  
 \pf    The horizontal composition specified by (\ref{E:horcompMorGV}) is a genuine composition law for morphisms in the sense that it is the composition law in the product category $\mbg\times\mbV$, where source and target maps given by
  $$s_h(k,f)=\bigl(s(k), s(f)\bigr)\qquad\hbox{and}\qquad t_h(k,f)=\bigl(t(k), t(f)\bigr).$$
 
It is readily checked that the composition  (\ref{E:compMorGV}) has the correct behavior with regard to source and target maps, and is associative. Furthermore, for each $f\in\Mor(\mbV)$, the element $(1_e, f)$ is the identity morphism $f\to f$, where $1_e\in\Mor(\mbG)$ is the identity element in the morphism group.  

It remains only to verify the exchange law. On one hand we have
\begin{eqnarray}\label{E:exchVG}
\bigl((k'_2,f'_2)\circ (k'_1, f'_1)\bigr)\circ_h \bigl((k_2, f_2)\circ (k_1, f_1)\bigr) & &\\
&&\hskip -1in= (k'_2k'_1, f'_1)\circ_h (k_2k_1, f_1)  \nonumber\\
&&\hskip -1in=\bigl((k'_2k'_1)\circ (k_2k_1), f'_1\circ f_1\bigr)  \nonumber\\
&&\hskip -1in=\bigl((k'_2\circ k_2)(k'_1\circ k_1), f'_1\circ f_1\bigr)  \nonumber\\
&&\hbox{(using the exchange law (\ref{E:exchange}))}, \nonumber
\end{eqnarray}
and on the other
\begin{eqnarray}
\bigl((k'_2,f'_2)\circ_h(k_2,f_2)\bigr)\circ \bigl((k'_1,f'_1)\circ_h(k_1,f_1)\bigr)&\nonumber\\
&&\hskip -1in= (k'_2\circ k_2, f'_2\circ f_2)\circ (k'_1\circ k_1, f'_1\circ f_1) \nonumber\\
&&\hskip -1in =\bigl((k'_2\circ k_2)(k'_1\circ k_1), f'_1\circ f_1\bigr),\nonumber
\end{eqnarray}
in agreement with the last line in (\ref{E:exchVG}).
 \epf 

\section{Crossed modules and representations}\label{s:matrep}

  In this section we   will use the   crossed module  $(G, H, \alpha,\tau)$ structure for a categorical group $\mbg$  (as discussed for  (\ref{E:catcross})) to   construct an action of a categorical group $\mbG$ from a representation of the underlying object group $G=\Obj(\mbG)$. 
 
Let $V$ be a vector space and ${\rm End}(V)$ be the set of endomorphisms of $V$. Let us define a category  $\mbV$, whose
object set is $V$ and   a morphism  is given by an ordered pair $(f, v)$, where $f\in {\rm End}(V)$ and $v\in V$. Source and targets are given by
\Beq \label{E:autact}
s(f,v)=v\qquad\hbox{and}\qquad t(f,v)=f(v),
\Eeq
and   composition of morphisms is given by
\Beq \label{E:autcom}
  (f',v')\circ (f,v)= (f'\circ f,v),
\Eeq
defined when $v'=f(v)$. 

  Let   ${\mathcal D}$ be a representation of $G$ on a vector space $V$:
$${\mathcal D}:G\to \End(V): g\mapsto {\mathcal D}_g.$$ 
Composing with the  homomorphism $\tau:H\longrightarrow G$  gives a representation
$T$ of $H$ on $V$, given by
\Beq \label{E:repT}
T_h=\mcD_{\tau(h)},
\Eeq
for every $h\in H$. Now we construct a functor
$$\rho: \mbG\times\mbV\to\mbV.$$
On objects it is given by
$$\rho(g)u= \mclD_g(u).$$
On morphisms it is given by
\begin{eqnarray} \label{E:endmap}
\rho  :\Mor\bigl(\mbg\bigr)\times \Mor\bigl(\mbV\bigr) &\longrightarrow & \Mor\bigl(\mbV\bigr)\\
 \bigl((h,g),\,(u, f)\bigr) &\mapsto & \bigl(T_h\mcD_g f {{\mcD}_g^{-1}}, \, \mcD_g(u)\bigr) \nonumber\\
 &&=\bigl(\mcD_{g'}f\mcD_g^{-1},\,\mcD_gu\bigr), \nonumber
\end{eqnarray}
for all  $g\in G, h\in H$, $u\in V$ and $f\in {\rm End}(V)$, with $g'=\tau(h)g=t(h,g)$; we have used
\begin{equation}\label{E:Thgpg}
T_h \mcD_g =\mcD_{g'}.
\end{equation} 
Diagrammatically put, the action of $\rho(h,g)$ on $(u,f)$ may be described by

\begin{equation} \label{D:rhohguf}
 \xymatrix{
        u_1 \ar@{-}[r] & f\ar[d]_{\rho(h,g)}\ar[r] & u_2         \\
   \mcD_{g}u_1 \ar@{-}[r] & \mcD_{g'}f\mcD_g^{-1}\ar[r] & \mcD_{g'}u_2     }  
   \end{equation}

 We show now that  $\rho$ is indeed a functor and provides an action of $\mbG$.

\begin{prop}\label{pr:repsimlcm}
Let $\mbg$ be the categorical group associated with a crossed module ${\ghat}$. Let $\mbV$ be the   category with $\Obj\bigl( \mbV \bigr)=V$, a vector space, and for which a morphism $u_1\to u_2$ is 
given by $(u_1, f)$, where $f\in {\rm End}(V)$ and $u_2=f(u_1)$. The composition 
 is as in (\ref{E:autcom}). Let 
$\mcD$ be a representation of $G$ on $V$. Then $(g, u)\mapsto \mcD_gu$
and $\rho$ as defined in (\ref{E:endmap}) give an action of  $\mbg$ on $\mbV$.
\end{prop}
\pf  For notational simplicity we shall write $a\cdot x$ for $\rho(a)x$, where $a$ and $x$ may be objects or morphisms. Moreover, composition of endomorphisms of $V$ will also be written as juxtaposition, as already done in (\ref{E:endmap}). The arguments below can be understood by consulting the diagram:

\begin{equation}
    \xymatrix{
        u_1 \ar@{-}[r] & {f_1}\ar[d]_{(h_1,g_1)}\ar[r] & u_2  \ar@{-}[r]&  f_2\ar[d]_{(h_2,g_2)} \ar[r]& u_3 \\
  \mcD_{g_1}u_1 \ar@{-}[r] & {\mcD_{g_2}f_1\mcD_{g_1}^{-1}}\ar[r] & \mcD_{g_2}u_2  \ar@{-}[r]&  \mcD_{g_3}f_2\mcD_{g_2}^{-1} \ar[r]& \mcD_{g_3}u_3 \\     }
\end{equation}

Using  (\ref{E:crossedprod}),  (\ref{E:crossedcomp}) and (\ref{E:autcom}) it is a straightforward verification that, whenever well defined,
\begin{eqnarray*}
 \bigl((h_2, g_2)\circ (h_1,g_1)\bigr)\cdot\bigl(  (f_2, u_2)\circ (f_1, u_1)\bigr)&&\\
 &\hskip -1.25in=&\hskip -.5in
 (h_2h_1, g_1)\cdot (f_2\circ f_1, u_1) \quad\hbox{from  (\ref{E:crossedcomp}) and (\ref{E:autcom})}\\
&\hskip -1.25in=&\hskip -.5in
\bigl( T_{h_2h_1}\mcD_{g_1}f_2f_1\mcD_{g_1}^{-1}, \mcD_{g_1}(u_1) \bigr), \quad\hbox{by  (\ref{E:endmap}),} 
\end{eqnarray*}
which agrees with
\begin{eqnarray*}
\bigl((h_2,g_2)\cdot (f_2,u_2)\bigr)\circ \bigl((h_1,g_1)\cdot (u_1,f_1)\bigr) &\\
&\hskip -1.25in=&\hskip -.5in
\bigl(T_{h_2} \mcD_{g_2} f_2 {{\mcD}_{g_2}^{-1}}, \mcD_{g_2}u_2\bigr)\circ \bigl(T_{h_1} \mcD_{g_1} f_1 {{\mcD}_{g_1}^{-1}}, \mcD_{g_1}u_1 \bigr)\\
&\hskip -1.25in=&\hskip -.5in
\bigl(T_{h_2} \mcD_{g_2} f_2 {{\mcD}_{g_2}^{-1}}T_{h_1} \mcD_{g_1} f_1 {{\mcD}_{g_1}^{-1}}, \mcD_{g_1}u_1\bigr)\\
\end{eqnarray*}
because
\begin{eqnarray*}
T_{h_2} \mcD_{g_2} f_2 {{\mcD}_{g_2}^{-1}}T_{h_1} \mcD_{g_1} f_1 {{\mcD}_{g_1}^{-1}} &  =&T_{h_2} \mcD_{g_2}f_2\mcD_{g_2}^{-1}\mcD_{\tau(h_1)g_1}f_1\mcD_{g_1}^{-1}\\
&=&T_{h_2} \mcD_{\tau(h_1)g_1}f_2f_1\mcD_{g_1}^{-1}\\
& =& T_{h_2h_1}\mcD_{g_1}f_2f_1\mcD_{g_1}^{-1}.
\end{eqnarray*}
Hence $\rho$ is a functor. Next, we check that $\rho$ is a representation. Consider $(h_1,g_1), (h_2, g_2)\in H\times G$ that are composable as elements of $\Mor(\mbG)$; this means
$$g_2=\tau(h_1)g_1.$$

Then

\begin{eqnarray*}
\bigl((h_2,g_2) (h_1,g_1)\bigr)\cdot (f,u) &=&\Bigl(h_2\alpha(g_2)(h_1), g_2g_1\Bigr)\cdot (f,u)\\
& =&\bigl(T_{h_2\alpha(g_2)(h_1)}\mcD_{g_2g_1}f  \mcD_{g_2g_1}^{-1}, \,\mcD_{g_2g_1}u\bigr)
\end{eqnarray*}

which agrees with

\begin{eqnarray*}
 (h_2,g_2)\cdot\bigl[(h_1,g_1)\cdot (f,u)\bigr] &=&(h_2,g_2)\cdot\bigl(T_{h_1}\mcD_{g_1}f\mcD_{g_1}^{-1},\, \mcD_{g_1}u\bigr)\\
&= &\bigl( T_{h_2}\mcD_{g_2} T_{h_1}\mcD_{g_1}f\mcD_{g_1}^{-1}\mcD_{g_2}^{-1}, \,\mcD_{g_2}\mcD_{g_1}u\bigr),\\
\end{eqnarray*}
on using the Peiffer  identity $\tau(\alpha(g_2)(h_1))=g_2\tau(h_1)g_2^{-1}$ noted earlier in  (\ref{E:Peiffer}).  \epf

\section{Twists and semidirect products of categories}\label{s:semidir}

In this section we introduce   the notion of   {\em semidirect product} or {\em $\eta$-twisted product} of
categories.  It should be noted here that the namesake of semidirect product of categories was
introduced in \cite{steinberg}. However our definition of semidirect product is different both in approach and motivation. Our primary motivation comes from the discussion in the Introduction in the context of (\ref{E:compsemid}).

Let   ${\mathbf A}$ and ${\mathbf B}$ be categories. Let  
$$\eta:\Mor\bigl({\mathbf  A}\bigr)\times  \Mor\bigl({\mathbf B}\bigr)\longrightarrow \Mor\bigl({\mathbf A}\bigr)$$
be a map satisfying the following conditions:
\begin{enumerate}
\item For any    $k\in \Hom(a_1,a_2)$ and $f\in \Hom(b_1,b_2)$ ,
\Beq \label{etacon1}
\eta(k,f)\in \Hom(a_1, a_2),
\Eeq
  where $a_1,a_2\in \Obj\bigl({\mathbf A}\bigr), b_1,b_2\in \Obj\bigl({\mathbf B}\bigr).$
\item  The map $\eta$ behaves as a right action: for any $k\in \Mor\bigl({\mathbf A}\bigr)$ and 
 composable $f_2, f_1\in \Mor\bigl(\mathbf B\bigr)$:

\begin{eqnarray} \label{etacon2}
\eta\bigl(k,f_2\circ f_1  \bigr) &=&\eta\bigl(\eta(k,f_2),f_1 \bigr) \\
\eta(k,1_b ) &=&k \nonumber
\end{eqnarray}
where $1_b$ is the identity morphism $b\to b$ for any $b\in \Obj(\mbbb)$.

 \item For any $f\in \Mor\bigl({\mathbf B}\bigr)$ and composable $k_1,k_2\in \Mor\bigl({\mathbf A}\bigr)$
\begin{eqnarray}\label{etacon3} 
\eta \bigl(k_2\circ k_1,f\bigr)&=&\eta \bigl(k_2,f\bigr)\circ\eta\bigl(k_1,f\bigr)\\
\eta(1_a,f) &=& 1_a\nonumber
\end{eqnarray}
for every object $a\in\Obj({\mathbf A})$.

\end{enumerate}

 
\begin{prop}\label{pr:semidir}
Let ${\mathbf A}$ and ${\mathbf B}$ be categories and $\eta:\Mor\bigl({\mathbf A}\bigr)\times \Mor\bigl({\mathbf B}\bigr)\longrightarrow \Mor\bigl({\mathbf A}\bigr)$ satisfy
conditions ${\rm {1-3}}$ listed above. Then $\eta$ specifies a category ${\mathbf A}\rtimes_{\eta} {\mathbf B}$ that has the following description:
\begin{eqnarray}
&&\Obj({\mathbf A}\rtimes_{\eta} {\mathbf B}):=\Obj({\mathbf A}) \times \Obj({\mathbf B})\label{objsemi}\\
&&\Mor({\mathbf A}\rtimes_{\eta} {\mathbf B}):=\Mor({\mathbf A}) \times\Mor({\mathbf B}),\label{morsemi}
\end{eqnarray}
target and source maps are the obvious ones, and the composition law is given by
\begin{equation}\label{E:semicompo} 
(k_2,f_2)\circ_{ {\eta}}(k_1,f_1):=\bigl(\eta(k_2,f_1)\circ k_1,f_2\circ f_1\bigr). 
\end{equation}
\end{prop}
 We call ${\mathbf A}\rtimes_{\eta} {\mathbf B}$ the {\em semidirect product}  of the categories ${\mathbf A}$ and ${\mathbf B}$
with respect to $\eta$. In the special case $\eta:(k,f)\mapsto k$ the semi-direct product reduces to the ordinary (direct) product.

 \pf  
 Condition (\ref{etacon1}) ensures that composition in the right hand side of (\ref{E:semicompo}) is well defined. To verify associativity, consider composable morphisms  $(k_3, f_3)$, $(k_2, f_2)$ and $(k_1, f_1)$  in ${\mathbf A}\rtimes_{\eta} {\mathbf B}$; then  
\begin{eqnarray}\label{asso1}
 \bigl((k_3,f_3)\circ_{ {\eta}}(k_2, f_2)\bigr)\circ_{ {\eta}}(k_1, f_1)& =& \bigl( \eta(k_3,f_2)\circ k_2,  f_3\circ f_2\bigr) \circ_{ {\eta}}(k_1,f_1)\\
&=& \left(  \eta\bigl(\eta(k_3,f_2)\circ k_2,f_1\bigr)  \circ k_1, f_3\circ f_2\circ f_1\right)\nonumber
\end{eqnarray}
and 
\begin{eqnarray}\label{asso2}
 (k_3, f_3)\circ_{ {\eta}}\Bigl((k_2, f_2)\circ_{ {\eta}}(k_1, f_1)\Bigr)   & =&  (k_3, f_3)\circ_{ {\eta}} \Bigl( \eta(k_2,f_1)\circ k_1,f_2\circ f_1\Bigr) \nonumber\\
& =& \Bigl(\eta\bigl(k_3,f_2\circ f_1\bigr)\circ \eta(k_2,f_1)\circ k_1, f_3\circ f_2\circ f_1,  \Bigr).\nonumber
\end{eqnarray}
 
  Conditions (\ref{etacon2}) and  (\ref{etacon3}) imply  that the right hand sides
of (\ref{asso1}) and (\ref{asso2}) are equal. These conditions also clearly imply that $(1_a, 1_b):(a,b)\to (a,b)$ is the identity morphism associated with 
$(a,b)\in \Obj\bigl({\mathbf A}\times {\mathbf B}\bigr),$ where $a\in \Obj\bigl({\mathbf A}\bigr)$ and $b\in \Obj\bigl({\mathbf B}\bigr).$\epf

{\bf Example SD1.} A homomorphism $\alpha:G_2\longrightarrow {\rm Aut}(G_1):b\mapsto\alpha_b$, where $G_1$ and $G_2$ are groups, specifies a semidirect
product group  $G_1\rtimes_\alpha G_2$ given by the group product
$$(a_2, b_2)(a_1, b_1)=\bigl(a_2\alpha_{b_2}(a_1), b_2 b_1\bigr),$$
for any $a_1, a_2 \in G_1$, $b_1, b_2 \in G_2$. Let  $\mbg_1$ be the usual   single-object  category 
associated with the group $G_1$, and $\mbg_2$ be   likewise
for $G_2$; in particular, $\Mor(\mbg_j)=G_j$. Then  $G_1\rtimes_{\alpha}G_2$ naturally defines the semidirect product category
$\mbg_1\rtimes_{\eta_{\alpha}}\mbg_2$, where
$$\eta_{\alpha}: G_1\times G_2\to G_1: (a,b)\mapsto \alpha_b(a).$$
 Conversely, given a semidirect product between two single 
object categories we recover the semidirect product of groups. ${ }$

{\bf Example SD2.} Let us again consider the crossed module ${\ghat}$ and the corresponding categorical group $\mbg$. 
Suppose  morphisms $(h,g), (h',g)$ have the same source and target elements. Thus $\tau(h)g=\tau(h')g$; hence
$h$ and $h'$ must be related as $h'=a h $, for some $a \in \ker {(\tau)}$. Now let $V$ be a vector space. Define
a single object category ${\mathbf C}({V})$, the single object being $V$ and a morphism 
of  ${\mathbf C}({V})$ being any element of ${\rm Aut}(V)$, the set of all invertible endomorphisms of $V$; thus: 
$$\Obj\bigl({\mathbf C}({V})\bigr)=\{V\}, \quad\Mor\bigl({\mathbf C}({V})\bigr)={\rm Aut}(V).$$
Let $\mu:{\rm Aut}(V)\longrightarrow {\rm Ker}(\tau)\subset H$ be a homomorphism.  Then we obtain a homomorphism 
${\rm Aut}(V)\longrightarrow {\rm Aut}(H)$ given by conjugation $$\bigl(\mu(f)\bigr)(h):=\mu(f)^{-1} h \mu(f):=\mu(h,f),$$
 where $\mu(f)\in {\rm Ker}(\tau)$.
 We will now define a semidirect product between the categories $\mbg$ and ${\mathbf C}{(V)}$. Define a map 
$$\eta: \Mor\bigl(\mbg\bigr)  \times\Mor\bigl({\mathbf C}({V})\bigr) \longrightarrow \Mor\bigl(\mbg\bigr)$$
given by 
$$\bigl((h,g); f\bigr)\mapsto \bigl( \mu(h,f),g\bigr)\in H\rtimes_{\alpha}G\simeq \Mor(\mbG),$$
where $g\in G, h\in H$ and $f\in \Mor\bigl({\mathbf C}({V})\bigr)={\rm Aut}(V)$. Both $\bigl(\mu(h,f),g\bigr)$ and $(h,g)$ have source $g$; they also have a common target:
$$\tau\left(\mu(h,f)\right)g=\tau\Bigl(\mu(f)^{-1}\Bigr)\tau(h)\tau\Bigl(\mu(f)\Bigr)g=\tau(h)g$$ 
because $\mu(f)$ is assumed to be in $\ker\tau$.  Next  we check that
$$\eta\bigl((h,g); f_2\circ f_1 \bigr)= \bigl(\mu(f_2f_1)^{-1}h\mu(f_2f_1),g\bigr)$$
agrees with
\begin{eqnarray}
\eta\Bigl( \eta\bigl((h,g); \, f_2\bigr), f_1\Bigr) &= &\eta\Bigl( \bigl(\mu(f_2)^{-1}h\mu(f_2),g\bigr);\, f_1\Bigr)\\
&=& \Bigl(\mu(f_1)^{-1}\mu(f_2)^{-1}h\mu(f_2)\mu(f_1),\,g\Bigr),\nonumber
\end{eqnarray}
and, with $g'=\tau(h)g$, we have

\begin{eqnarray}
\eta\bigl((h',g')\circ (h,g);\, f \bigr) &=&\eta\Bigl((h'h,g); \,f\Bigr)\qquad\hbox{(using (\ref{E:crossedcomp}))}
\\
&=&\Bigl(\mu(f)^{-1}h'h\mu(f),g\Bigr) \nonumber\\
&=&\eta\bigl( (h',g'); f \bigr) \circ \eta\bigl( (h,g); f \bigr).\nonumber \end{eqnarray}
 
Hence $\eta$ satisfies conditions (\ref{etacon1})--(\ref{etacon3}). By Proposition~\ref{pr:semidir} it
defines a semidirect product category $ \mbg\rtimes_{\eta}{\mathbf C}({V}) $. The composition, stated explicitly, is
\begin{eqnarray} \label{ex:semi2} \bigl((h_2, g_2); f_2\bigr)\circ_{ {\eta}}\bigl((h_1, g_1); f_1\bigr) &=&\bigl(\eta\bigl((h_2, g_2), f_1\bigr)\circ (h_1, g_1); f_2\circ f_1\bigr)\\
&=& \Bigl(\bigl(\mu(f_1)^{-1}h_2\mu(f_1)h_1, g_1\bigr); f_2\circ f_1\Bigr).  {\fbox{ }}\nonumber \end{eqnarray}

 {\bf Example VGG}.
 Let $\alpha:G\to \End(V)$ be a representation of a  group $G$ on a vector space $V$ over some field. Then, with $V$ viewed as an additive group, we have  a crossed module
 $$(G, V, \alpha,\tau),$$
 where the target map $\tau$ is the constant map $\tau:V\to G: v\mapsto e$. 
 We  have discussed this in Example CG2 in section \ref{S:catgrpsumm}. The object group of this categorical group is $G$ and the morphism group is the semi-direct product $V\rtimes_{\alpha} G$.  There is another example we can construct out of this.  Consider the category $\mbV$, whose object set is the one-element set $\{V\}$ and whose morphism set is $V$, with composition being given by addition of vectors. Next consider the categorical group $\mbg_0$, whose object group is $G$ and a unique morphism $(b,a)$   goes from each (source) $a\in G$ to each (target) $b\in G$.  Consider next the map
 $$\eta: V\times (G\times G)\to V:\Bigl(v, (g_2,g_1)\Bigr)\mapsto \alpha(g_1g_2^{-1})v=g_1g_2^{-1}v.$$
 We write $\alpha(g)w$ as $gw$, and verify that
 \begin{eqnarray*}
 \eta\Bigl(\eta\bigl(v, (g_3,g_2)\bigr),   (g_2,g_1)\Bigr) &= &\eta\Bigl(g_2g_3^{-1}v, (g_2,g_1)\Bigr)\\
 &=&g_1g_2^{-1}g_2g_3^{-1}v=g_1g_3^{-1}v\\
 &=&\eta\Bigl(v, (g_3,g_2)\circ (g_2,g_1)\Bigr) 
 \end{eqnarray*}
 and
 \begin{eqnarray*}
  \eta\bigl(v', (g_2,g_1)\bigr) \circ\eta\Bigl(v,   (g_2,g_1)\Bigr) &=& g_2g_1^{-1}v'+g_2g_1^{-1}v\\
  &=&g_2g_1^{-1}(v'+v)\\
  &=&\eta\bigl(v'\circ v, (g_2,g_1)\bigr),
  \end{eqnarray*}
for all $g_1,g_2,g_3\in G$ and all $v,v'\in V$.
Moreover,
\begin{eqnarray*}
\eta\bigl(v, (g,g)\bigr) &=& v\\
\eta\bigl(0,(g_2,g_1)\bigr) &=&0
\end{eqnarray*}
for all $v\in V$ and all $g,g_1,g_2\in G$. Thus $\eta$ satisfies  the conditions (\ref{etacon1}), (\ref{etacon2}),(\ref{etacon3}). The composition law of morphisms for the semidirect product category $\mbV\rtimes_{\eta}\mbg_0$ is

\begin{eqnarray}\label{E:compsemiGV}
\Bigl(w, (g_3,g_2)\Bigr)\circ \Bigl(v, (g_2,g_1)\Bigr) &=&\Bigl(\eta\bigl(w, (g_2,g_1)\bigr)+v,  (g_3,g_1)\Bigr)\\
&=&
\Bigl(g_1g_2^{-1}w+v, (g_3,g_1)\Bigr).\nonumber \hskip 1in {\fbox{ }}
\end{eqnarray}

{\bf Example 1D}. Let $G$ be a group, $\alpha:G\to\End(H)$  a representation on a finite dimensional complex vector space $H$, and $\lambda_0$  a one-dimensional representation of $H$; we write $\alpha(g)h$ simply as $gh$. Let $\mbg$ be the categorical group whose object group is $G$, the morphism group is $H\times_{\alpha} G$, source and targets are given by
$$s(h,g)=t(h,g)=g$$
and composition by
\begin{equation}\label{E:comph1h2g}
(h_2,g)\circ (h_1,g)=(h_1+h_2, g).
\end{equation}
Consider
\begin{equation}\label{E:eta0GH}
\eta_0: V\times (H\rtimes_{\alpha} G)\to V: (v; h,g)\mapsto \lambda_0(g^{-1}h) v.
\end{equation}
We note that $\eta_0(0;h,g)=0$ and $\eta_0(v;0,e)=v$.
Working with $g_2=g_1$,
\begin{eqnarray*}
\eta_0\bigl(\eta_0(v;h_2,g_2); (h_1,g_1)\bigr) &=&\eta_0\bigl(\lambda_0(g_2^{-1}h_2)v; h_1, g_1\bigr)\\
&=&\lambda_0(g_1^{-1}h_1)\lambda_0(g_2^{-1}h_2)v\\
&=&\lambda_0\bigl(g_1^{-1}(h_1+h_2)\bigr)v\\
&&\quad\qquad \hbox{(using $g_1=g_2$)}\\
&=&\eta_0\bigl(v; (h_2,g_2)\circ (h_1,g_1)\bigr)
\end{eqnarray*}
and

\begin{eqnarray*}
 \eta_0(v_2;h,g)\circ \eta_0\bigl(v_1; h,g\bigr)  &=&\lambda_0(g^{-1}h)v_2+\lambda_0(g^{-1}h)v_1\\
&=&\lambda_0(g^{-1}h)(v_2+v_1)\\
&=& \eta_0(v_2\circ v_1;h,g).
\end{eqnarray*}
Thus $\eta_0$ specifies a semidirect product $\mbV\rtimes_{\eta_0}\mbg$, where $\mbV$ is the category with one object and with morphism group being the additive group of $V$. ${\fbox{ }}$

It has been pointed out to us by the referee that in the preceding example $\lambda_0$ need not be $1$-dimensional. 

\section{Semidirect products and  twisted actions}\label{s:semirepre}

In this section we develop a notion of  a `twisted action' of a categorical group on a category, the twist arising from an $\eta$-map as discussed in the preceding section.

Let ${\mathbf V}$ be a   category, $\mbg$   a categorical group, and suppose 
$$\eta: \Mor\bigl(\mbg \bigr)\times \Mor\bigl(\mbV \bigr)   \longrightarrow \Mor\bigl(\mbg \bigr)$$
satisfies conditions (\ref{etacon1})--(\ref{etacon3}). By Proposition~\ref{pr:semidir}
there is a semidirect product category $\mbg \rtimes_{\eta} \mbV$.  The composition law for morphisms is
\begin{equation}\label{E:semidirleftcompo}
(k_2,f_2)\circ_ {\eta}(k_1,f_1)=\left( \eta(k_2,f_1)\circ k_1,f_2\circ f_1\right)
\end{equation}
whenever the right side is well-defined. In this section we will assume, furthermore, that
\begin{eqnarray}\label{E:etak1g} 
\eta(k1_a,f) &=&\eta(k,f)1_a\\
\eta(1_ck,f) &=&1_c\eta(k,f) \nonumber
\end{eqnarray}
for all $k\in\Mor(\mbg)$, $a,c\in \Obj(\mbg)$, and $f\in\Mor(\mbV)$. Then for any 
$k_1:a\to b$ and $k_2:c\to d$ in $\Mor(\mbg)$ and any $f\in\Mor(\mbV)$,
\begin{eqnarray}\label{E:etaprodk2k1}
\eta(k_2k_1,f)&=&\eta\left( (k_2\circ 1_c)(1_b\circ k_1), f\right)\\
&=&\eta\left((k_21_b)\circ (1_ck_1), f\right)\nonumber\\
&=&\eta(k_21_b,f)\circ\eta(1_ck_1, f)\qquad\hbox{(using (\ref{etacon3}))}\nonumber\\
&=& [\eta(k_2,f)1_b]\circ [1_c\eta(k_1,f)]\nonumber\\
&=& [\eta(k_2,f)\circ 1_c][1_b\circ \eta(k_1,f)]\nonumber\\
&=&\eta(k_2,f)\eta(k_1,f).\nonumber
\end{eqnarray}
Thus $\eta$ respects not only the composition law in $\Mor(\mbG)$ but also the group product operation.

By an {\em action} of $\mbg$ associated to $\eta$, or an {\em $\eta$-twisted action} of $\mbg$, we mean a functor 
\begin{equation} \mbrho:\mbg \rtimes_{\eta} \mbV \longrightarrow \mbV  \label{def:semirepfunc}
\end{equation}
for which
\begin{equation}\label{E:rhosd1ef}
\rho(1_e,f)=f
\end{equation}
for all $f\in\Mor(\mbV)$, with $e$ being the identity element in $\Obj(\mbg)$, and
\begin{equation}\label{E:rhosdk2k1f}
\rho\left(k_2,\rho(k_1,f)\right)=\rho(k_2k_1,f)
\end{equation}
for all $k_1, k_2\in\Mor(\mbg)$ and $f\in\Mor(\mbV)$. As seen earlier in the context of (\ref{E:rhoevv}) and (\ref{E:rhoabv}) we then have
\begin{eqnarray*}
\rho(e,v) &=&v\\
\rho\left(a,\rho(b,v)\right)&=&\rho(ab,v)
\end{eqnarray*}
for all $a,b\in\Obj(\mbg)$ and all $v\in \Mor(\mbV)$.

The {\em key new idea built into this notion is encoded in the requirement  that $\rho$ be   also  a functor}; this  implies that 
\begin{equation}\label{E:rhoetak2f1} 
 \rho(k_2,f_2)\circ\rho(k_1,f_1)=\rho\left(\eta(k_2,f_1)\circ k_1, f_2\circ f_1\right)   
\end{equation}
for all   $k_1,k_2\in\Mor(\mbG)$, and $f_1,f_2\in\Mor(\mbV)$ for which the composites $k_2\circ k_1$ and $f_2\circ f_1$ are defined.  Dropping the  explicit reference to $\rho$, this becomes:
\begin{equation}\label{E:rhoetak2f1b} 
(k_2\cdot f_2)\circ (k_1\cdot f_1)=\bigl(\eta(k_2,f_1)\circ k_1\bigr)\cdot\bigl(f_2\circ f_1\bigr),   
\end{equation}
where $a\cdot x$ means $\rho(a,x)$.

For the special case
$\eta={\rm Pr_1}:(k, \phi)\mapsto k$   we recover   actions as defined in (\ref{def:repfunc}.  The distinction between this simplest case and the general semidirect product is  expressed in the  functoriality of $\rho$, as explicitly given through (\ref{E:rhoetak2f1}).

The following is the analog of Proposition \ref{prop:doublerep}:

\begin{prop}\label{prop:doublesemirep}
Suppose $\rho:\mbg\rtimes_{\eta}\mbV\to\mbV$ is an action associated to $\eta$, where $\mbg$ is a categorical group, $\mbV$ is a  category,   $\eta$   satisfies the relations (\ref{E:etak1g}) as well as
\begin{equation}\label{E:etalkf}
\eta(k,k'f)=\eta(k,f)
\end{equation}
for all $k,k'\in\Mor(\mbg)$ and all $f\in\Mor(\mbV)$.  Then there is a double category over $\mbV$ whose arrows are the elements $(k,f)\in \Mor(\mbg)\times\Mor(\mbV)$, with source and target  given by
\begin{equation}\label{E:Pdoublesemi}
s(k,f)=f\,\quad \hbox{and}\quad t(k,f)=\rho(k,f),
\end{equation}
composition  given by
\begin{equation}
(k',f')\circ (k,f)=(k'  k, f),
\end{equation}
whenever $f'=\rho(k')f$, and horizontal composition  given by the semidirect product composition
\begin{equation}\label{E:Pdoublesemicomp}
(k_2,f_2)\circ_h (k_1,f_1)=(k_2,f_2)\circ_ {\eta} (k_1,f_1)
\end{equation}
whenever the composition on the right is defined.
\end{prop}

We note that $\circ_h$ is a genuine composition law in a category, and source and target maps in this category are given by
$$s_h(k,f)=\bigl(s(k), s(f)\bigr)\qquad\hbox{and}\qquad t_h(k,f)=\bigl(t(k), t(f)\bigr),$$
as in (\ref{prop:doublerep}). 
 We shall denote the double category specified in this Proposition by 
$$\Delta^\eta_{\rho}(\mbV).$$
The relation (\ref{E:etalkf}) expresses an invariance of $\eta$ with respect to $\rho$.

\pf First let us check that the source and  target of horizontal composition $\circ_h$ behave correctly. We have 
 \begin{eqnarray*}
 s_h[(k_2,f_2)\circ_h (k_1,f_1)]&=&s_h\left(\eta(k_2,f_1)\circ k_1, f_2\circ f_1\right)\\
 &= &\bigl(s(k_1),   s( f_1)\bigr)\\
 &=& s_h(k_1,f_1)
 \end{eqnarray*}
 and
  \begin{eqnarray*}
 t_h\bigl[(k_2,f_2)\circ_h (k_1,f_1)\bigr] &=&t_h\left(\eta(k_2,f_1)\circ k_1, f_2\circ f_1\right)\\
 &=& \Bigl(t\bigl(\eta(k_2,f_1\bigr)\circ k_1\bigr), t(f_2\circ f_1)\Bigr)\\
 &=&\bigl(t(k_2), t(f_2)\bigr)\\
 &=& t_h(k_2, f_2).
 \end{eqnarray*}

The source and target maps for the composition law $\circ$ also behave correctly as we now verify. We have 
  \begin{eqnarray*}
 s[(k',f') \circ (k,f)]&= &s(k'k,f)\\
 & =& f\\
 &=& s(k,f)
 \end{eqnarray*}
 and
   \begin{eqnarray*}
 t[(k',f') \circ (k,f))]&=&\rho(k'k)f\\
 &=&\rho(k')\bigl(\rho(k)f\bigr)\\
 &=&\rho(k')f'\\
 &=&t(k',f').  
 \end{eqnarray*}

 We need now only verify  the exchange law. On one hand we have 
  \begin{eqnarray}
\label{E:exchVG2} 
\left[ (k'_2, f'_2)\circ (k_2,f_2)\right]\circ_h \left[(k'_1, f'_1)\circ (k_1, f_1) \right] & &\\
\hskip -1in&\hskip -1.5in  =&\hskip -.75in (k'_2k_2, f_2)\circ_h (k'_1k_1, f_1)\nonumber\\
\hskip -1in&\hskip -1.5in  =&\hskip -.75in\left(\eta(k'_2k_2, f_1)\circ(k'_1k_1), f_2\circ f_1\right)\nonumber 
\end{eqnarray}
and on the other 
\begin{eqnarray*}
\bigl((k'_2,f'_2)\circ_h(k'_1,f'_1)\bigr)\circ \bigl((k_2,f_2)\circ_h(k_1,f_1)\bigr)
& &\\
\hskip -1in&\hskip -1.5in  =&\hskip -.75in\bigl( \eta(k'_2,f'_1)\circ k'_1, f'_2\circ f'_1  \bigr) \, \circ \,\bigl( \eta(k_2,f_1)\circ k_1, f_2\circ f_1  \bigr) \nonumber\\
\hskip -1in&\hskip -1.5in  =&\hskip -.75in\Bigl( \bigl(\eta(k'_2,f'_1)\eta(k_2,f_1)\bigr)\,\circ (k'_1k_1), f_2\circ f_1\Bigr),\nonumber
\end{eqnarray*}
and so, noting that $f'_1=k_1f_1$ and  using the identities (\ref{E:etaprodk2k1})  and (\ref{E:etalkf}) satisfied by $\eta$ and comparing with (\ref{E:exchVG2}), we see that the exchange law holds.
  \epf

Let $(\eta_1, \rho_1, \mbV_1)$ and $(\eta_2, \rho_2, \mbV_2)$ be  actions of a categorical group $\mbg$ on    categories $\mbV_1$ and  $\mbV_2$, respectively. By a {\em morphism}
 $$\mbF: (\eta_1, \rho_1, \mbV_1)\longrightarrow (\eta_2, \rho_2, \mbV_2)$$
 we mean a functor 
$\mbF:\mbV_1 \longrightarrow \mbV_2$ that intertwines both   $\eta_j$ and $\rho_j$ in the sense that the following two diagrams commute:

\begin{equation} \label{D:etacom}
\xymatrix{
         \ar[d]^{{\rm Id}_{\mbg}\times {\mbF}} \Mor(\mbG)\times\Mor(\mbV_1)       \ar[r]^-{\eta_1} &\Mor(\mbG) \ar[d]_{{\rm Id}_{\mbg}} \\
 \Mor(\mbG)\times\Mor(\mbV_2) \ar[r]_-{\eta_2}&\Mor(\mbG)}
\end{equation}
and
 
\begin{equation} \label{D:semimap}
\xymatrix{
         \ar[d]^{{\rm Id}_{\mbg}\times {\mbF}} \mbG\rtimes_{\eta_1}\mbV_1    \ar[r]^-{\rho_1}&\mbV_1 \ar[d]_{ {\mbF}} \\
 \mbG\rtimes_{\eta_2}\mbV_2 \ar[r]_-{\rho_2}&\mbV_2 
 }
\end{equation}
 
 In (\ref{D:etacom}) arrows are set maps, and in (\ref{D:semimap}) arrows are functors.   
Thus the functor $\mbF$ must be `compatible' with maps $\eta_1$ and $\eta_2$.  In this context we make the following observation.

\begin{lemma}\label{l:compatietaf}
If diagram (\ref{D:etacom}) commutes then  ${\rm Id}_{\mbg}\times {\mbF}$
is a functor.
\end{lemma}
\pf
Let $k\in \Mor\bigl(\mbg\bigr)$ and $f\in \Mor\bigl(\mbV_1\bigr)$. Then the   commutation of the diagram (\ref{D:etacom}) implies
\Beq\label{firstcom}
\eta_2\bigl(k,  \mbF(f)\bigr)=\eta_1\bigl(k,  f \bigr).
\Eeq
At the morphism and object level
${\rm Id}_{\mbg}\times \mbF $ is defined as usual:
\Beqr
&&{\rm Id}_{\mbg}\times \mbF  :(k,f)\mapsto \bigl(k,\mbF(f)\bigr)\nonumber\\
&&{\rm Id}_{\mbg}\times \mbF:(g, v)\mapsto \bigl( g, \mbF(v)\bigr).\nonumber
\Eeqr
 For composable morphisms $(k', f')$ and $(k,f)$ in  ${\mbg\rtimes_{\eta_1} \mbV_1}$  we have, by the definition of the semidirect product structure in  (\ref{E:semicompo}),
$$
(k', f')\circ_{ {\eta_1}}(k, f):=\bigl(\eta_1(k',f)\circ k, f'\circ f \bigr).$$
Then we have
\begin{eqnarray*} 
({\rm Id}_{\mbg}\times \mbF)(k', f')\circ_ {\eta_2}({\rm Id}_{\mbg}\times \mbF)(k, f) & =&\left(\eta_2\bigl(k', \mbF( f)\bigr)\circ k, \mbF(f')\circ \mbF(f)   \right)\\
&=&\left(\eta_1(k', f)\circ k, \mbF(f'\circ f)  \right) \\
&&\qquad\hbox{ (using (\ref{firstcom}))}\\
&=&({\rm Id}_{\mbg}\times \mbF)\bigl(\eta_1(k',f)\circ k, f'\circ f\bigr)\\
&=&({\rm Id}_{\mbg}\times \mbF)\bigl((k', f')\circ_ {\eta_1}(k, f)\bigr),\end{eqnarray*} 
which shows that ${\rm Id}_{\mbg}\times\mbF$ is a functor.
\epf

Recall the invariance relation (\ref{E:etalkf}):
$$\eta_1\bigl(k, \rho_1(k',f_1)\bigr)=\eta_1(k, f_1),$$
for all $k, k'\in\Mor(\mbg)$ and all $f_1\in \Mor(\mbV_1)$,
that was needed in the double category result Proposition \ref{prop:doublesemirep}. We observe that this relation is preserved by morphisms; for, on using the commutativity in the diagrams (\ref{D:etacom}) and (\ref{D:semimap}), we have then
\begin{equation}\label{E:eta2Feta1}
\eta_2\bigl(k, \rho_2(k',f_2)\bigr)=\eta_2(k, f_2),
\end{equation}
for all $k, k'\in\Mor(\mbg)$ and all $f_2\in \Mor(\mbV_2)$ in the image of $\mbF$.

 Representations of $2$-groups on $2$-categories are studied in  \cite{BaezBFW}. Their definitions and development are different from ours. However, similar considerations arise there as well regarding  relationships between different representations. In their framework such relationships are expressed through $1$-intertwiners and their higher counterparts, called $2$-intertwiners. In our framework,   the corresponding notions are  functors between twisted actions and natural transformations between such functors.

Let $\mbF:(\eta_1, \rho_1, \mbV_1)\longrightarrow (\eta_2, \rho_2, \mbV_2)$ be a morphism between two actions 
of a categorical group $\mbg$.   Let $\Delta^{\eta}_{\rho}\bigl(\mbV\bigr)$ be the double category defined (in Proposition~\ref{prop:doublesemirep}) 
by the representation $(\eta, \rho)$ on $\mbV$. Then the morphism $\mbF:\bigl(\eta_1, \rho_1,\mbV_1\bigr) \longrightarrow \bigl(\eta_2, \rho_2,\mbV_2\bigr)$ 
induces a kind of functor between the double categories $\Delta^{\eta_1}_{\rho_1}\bigl(\mbV_1\bigr)\longrightarrow \Delta^{\eta_2}_{\rho_2}\bigl(\mbV_2\bigr)$. We explain this fully in:

\begin{prop}\label{pr:doubletodouble}
Let $\mbF:(\eta_1, \rho_1, \mbV_1)\longrightarrow (\eta_2, \rho_2, \mbV_2)$ be a morphism between actions  $(\eta_1, \rho_1)$ 
and $(\eta_2, \rho_2)$ of a categorical group $\mbg$ on categories $\mbV_1$ and $\mbV_2$ respectively. Let $\Delta^{\eta_1}_{\rho_1}\bigl(\mbV_1\bigr)$ and 
$\Delta^{\eta_2}_{\rho_2}\bigl(\mbV_2\bigr)$
 be the respective double categories over $\mbV_1$ and $\mbV_2$ defined 
by  $(\eta_1, \rho_1)$ and $(\eta_2, \rho_2)$ respectively. Then $\mbF$  induces a mapping $\mbF_*$, both on objects and on morphisms,
$$\Delta^{\eta_1}_{\rho_1}\bigl(\mbV_1\bigr)\longrightarrow \Delta^{\eta_2}_{\rho_2}\bigl(\mbV_2\bigr),$$
and $\mbF_*$ carries  vertical composition to vertical composition and horizontal composition to  horizontal composition.
\end{prop} 

\pf Recall that $\mbF$ is a functor
$$\mbF:\mbV_1\to\mbV_2.$$
  An object of  $\Delta^{\eta_1}_{\rho_1}\bigl(\mbV_1\bigr)$ is a morphism $f:v_1\to v'_1$ in $\Mor(\mbV_1)$. We define $\mbF_*$ on objects by
$$\mbF_*(f)=\mbF(f).$$
Next consider a morphism of $\Delta^{\eta_1}_{\rho_1}\bigl(\mbV_1\bigr)$. This is of the form
$$(k,f_1)\in \Obj(\mbg)\times\Mor(\mbV_1).$$
We define $\mbF_*$ on morphisms by
\begin{equation}\label{E:mbFstarkf1}
\mbF_*(k,f_1)=\bigl(k, \mbF(f_1)\bigr).
\end{equation}
Then for vertical composition we have
\begin{eqnarray*} 
\mbF_*\bigl((k',f_1')\circ (k,f_1)\bigr) &=&\mbF_*(k'k,f_1)\\
&=&\bigl(k'k, \mbF(f_1)\bigr)\\
&=&\mbF_*(k',f_1')\circ\mbF_*(k,f_1),
\end{eqnarray*}
where, for the last step, we note that the composition $\mbF_*(k',f_1')\circ\mbF_*(k,f_1)$ is meaningful because the target of $\mbF_*(k,f_1)$ is 
$$\rho_2\bigl(k,\mbF(f_1)\bigr),$$
which,  by the commutative diagram (\ref{D:semimap}), is equal to
$\mbF\bigl(\rho_1(k,f_1)\bigr)$, which is the same as $\mbF(f_1')$, the source of $\mbF_*(k',f_1')$.

Next we consider the effect of $\mbF_*$ on horizontal composition. Suppose
$(k,f)$ and $(k',f')$ are morphisms of $\Delta^{\eta_1}_{\rho_1}\bigl(\mbV_1\bigr)$ for which the composite
$(k',f')\circ_h(k,f)$ is defined; this means that the composite $f'\circ f$ is defined and also
$k'\circ k$ is defined. (Our goal at this point is to show first that the composite $\mbF(k',f')\circ_h\mbF(k,f)$ is defined.) Hence $\mbF_*(f')\circ\mbF_*(f)$ is defined and, moreover,
\begin{eqnarray}\label{E:Frho1kpfp} \mbF\bigl(\rho_1(k',f')\bigr)\circ\mbF\bigl(\rho_1(k,f)\bigr) &=& \rho_2\bigl(k',\mbF(f')\bigr)\circ \rho_2\bigl(k,\mbF(f)\bigr)\\
&=& \rho_2\bigl(\eta_2(k',\mbF(f))\circ k, \mbF\bigl(f'\circ f)\bigr)\nonumber
\end{eqnarray}
is also defined since $\eta_2(k',\mbF(f))$, having the same source and target as $k'$, is composable with $k$.  Thus, both the composites
$$\mbF(f')\circ\mbF(f)$$
and
$$\mbF\bigl(\rho_1(k',f')\bigr)\circ\mbF\bigl(\rho_1(k,f)\bigr)$$
are defined.
Hence the horizontal composite
$$\mbF_*(k',f')\circ_h\mbF_*(k,f) $$
is defined. Computing the effect of $\mbF_*$ on horizontal composition we have:
\begin{eqnarray}\label{E:Frho1kpfp2}
\mbF_*(k',f')\circ_h\mbF_*(k,f) &=&\bigl(k',\mbF(f')\bigr)\circ_h\bigl(k,\mbF(f)\bigr)\\
&=&\Bigl(\eta_2\bigl(k',\mbF(f)\bigr)\circ k, \mbF(f'\circ f)\Bigr) \nonumber\\
&=&\Bigl( \eta_1\bigl(k',f\bigr)\circ k, \mbF(f'\circ f)\Bigr)\qquad \hbox{(using (\ref{firstcom}))} \nonumber\\ 
&=&\mbF_*\bigl((k',f')\circ_h(k,f)\bigr).\nonumber
\end{eqnarray}
Thus $\mbF_*$ preserves both horizontal and vertical compositions.
\epf

Let $\mbF_1:(\eta_1, \rho_1, \mbV_1)\longrightarrow (\eta_2, \rho_2, \mbV_2)$ and $\mbF_2:(\eta_2, \rho_2, \mbV_2)\longrightarrow (\eta_3, \rho_3, \mbV_3)$
be two morphisms as described above. Then there is clearly a natural composition
$$\mbF_2\circ \mbF_1:(\eta_1, \rho_1, \mbV_1)\longrightarrow (\eta_3, \rho_3, \mbV_3),$$
and it is readily checked that this composition is associative. The identity functor ${\rm Id}_{\mbV}:\mbV\longrightarrow \mbV$ yields a morphism ${\rm Id}_\mbV:(\eta, \rho, \mbV)\longrightarrow (\eta, \rho, \mbV).$ 

\begin{prop}\label{pr:catrep}
There is a category ${\rm Rep}(\mbg)$ whose objects are twisted actions of a given categorical group $\mbg$ given by means of semidirect products with  categories and whose morphisms are the morphisms between such actions.   \end{prop}

 If $\mbF, \mbF', \mbF'':(\eta_1, \rho_1, \mbV_1)\longrightarrow (\eta_2, \rho_2, \mbV_2)$ and
$\omega:\mbF\Rightarrow \mbF'$, $\omega':\mbF'\Rightarrow \mbF''$ are natural transformations  then there is a composite natural transformation 
$\omega'\circ \omega:\mbF\Rightarrow \mbF''$ given by 
\begin{equation}\label{natver}
(\omega'\circ \omega) (v):=\omega'(v)\circ \omega(v) \qquad \forall v\in {\Obj}(\mbV_1) .
\end{equation}
Now suppose 
\begin{eqnarray}
&&\mbF_1, \mbF'_1:(\eta_1, \rho_1, \mbV_1)\longrightarrow (\eta_2, \rho_2, \mbV_2),\nonumber\\ 
&&\mbF_2, \mbF'_2:(\eta_2, \rho_2, \mbV_2)\longrightarrow (\eta_3, \rho_3, \mbV_3) \nonumber
\end{eqnarray}
are morphisms,
and $\omega_1:\mbF_1\Rightarrow \mbF'_1$, $\omega_2:\mbF_2\Rightarrow \mbF'_2$ are a pair of natural transformations; then   a natural transformation $\omega_2\circ_{H} \omega_1:\mbF_2\circ \mbF_1\Rightarrow \mbF'_2\circ \mbF'_1$ is specified by

\begin{eqnarray}\label{nathorz} 
(\omega_2\circ_{H} \omega_1)(v)&:= &{{\mbF}'_2}((\omega_1(v)))\circ \omega_2({\mbF}_1(v))\\
& &\rm {or \hskip .2 cm equivalently}\nonumber\\
&:= &\omega_2({\mbF}'_1(v))\circ {\mbF}_2((\omega_1(v))), \qquad \forall v\in {\Obj}(\mbV_1).\nonumber 
\end{eqnarray}

\section{An example with multiple twists}\label{S:P}

In this section we construct an example of twists and representations motivated by the Poincar\'e group $SO(1,3)\rtimes {\mathbb R}^4$. This group, being of great interest, has  been studied by others,  for example  \cite{Baez,BaezBFW,  BaezWise,  BarWise, CS, Pfeif}, in its role as a categorical group. Our study is entirely different in substance and flavor, and our objective is primarily to demonstrate how $\eta$-maps can be constructed in such natural but complex examples.

 In what follows one could take  $G$ to be the universal cover $SL(2,\mbc)$ of the proper orthochronous Lorentz group and  $H$ to be $\mbr^4$. The vector space $V$ used below can be thought of as a `spinorial' space on which a representation $S$ of $SL(2,\mbc)$ is given. 

Let $ \alpha:G\to {\rm End}(H)$ be a representation, where, as usual, $ G$ is the object group of a categorical  Lie group $ {\mathbf G}$, and $H$ is a real   vector space.  We write $\alpha$ as juxtaposition:
$$gh\stackrel{\rm def}{=}\alpha(g)h.$$
Then (as in Example CG2) we can form a categorical  Lie group 
$  {\mathbf G}$, associated with the Lie crossed module
$$(G, H, \alpha,\tau),$$
where $\tau$ is trivial; thus  the  object group is
$   {\rm Obj}({\mathbf G})=G $ 
and the morphism group is
$  {\rm Mor}({\mathbf G})=H\rtimes_{\alpha}G$. 
Composition of morphisms is given by (\ref{E:comph1h2g}); explicitly,
\begin{eqnarray}\label{E:compprodlaws} 
(h_2, g_2) (h_1, g_1) &=& (h_2+g_2h_1, g_2g_1)\\
(h_2, g)\circ (h_1, g) &=& (h_2 + h_1, g).\nonumber
\end{eqnarray}
Let $S:G\to\End(V)$ be a representation of $G$ on a   complex vector space $V$, and $M$ a manifold on which there is a smooth left action of $G$ (we can think of $M$ as either $\mbr^4$ or a hyperboloid of the form $\{p\in \mbr^4: \la p,p\ra_L=E^2\}$ for some fixed $E\in [0,\infty)$, with $\la\cdot,\cdot\ra_L$ being the Lorentz metric).  Now suppose  that  we have a character (one-dimensional representation) $\lambda_p$ of $H$, such that $(p,h)\mapsto\lambda_p(h)$ is smooth and
\begin{equation}\label{E:poinclambda}
\lambda_{gp}(h)=\lambda_p(g^{-1}h)
\end{equation}
for all $p\in M$, $h\in H$ and $g\in G$ (we think of the case $\lambda_p(h)=e^{i\la p,h\ra_L}$).  Let
$$\Psi$$
be the vector space of all smooth maps  $M\to V$ (in the physical context, such a map is  a matter field in momentum space).  
We define an action ${\tilde S}$ of $H\rtimes_{\alpha}G$ on $\Psi$ by:
\begin{equation}\label{E:poinc}
[(h,g)\cdot\psi](p)\stackrel{\rm def}{=}[{\tilde S}(h,g)\psi](p)\stackrel{\rm def}{=} \lambda_{p}(h)S(g)\psi(g^{-1}p) \end{equation}
for all $(h,g)\in H\rtimes_{\alpha}G$, $p\in M$, and $\psi\in \Psi$. We verify that this gives a representation: for any $(h_1, g_1), (h_2, g_2)\in H\rtimes_{\alpha}G$ and $p\in M$ we have
\begin{eqnarray}\label{E:poinc2} 
[(h_2,g_2)\bigl((h_1,g_1)\cdot\psi\bigr)](p) &=&\lambda_p(h_2)S(g_2)\bigl((h_1,g_1)\cdot\psi\bigr)](g_2^{-1}p)\\
&= &\lambda_p(h_2)S(g_2)S(g_1)\lambda_{g_2^{-1}p}(h_1)\psi(g_1^{-1}g_2^{-1}p)
\nonumber\\
&= &\lambda_p(h_2)\lambda_p(g_2h_1)S(g_2g_1)\psi\Bigl((g_2g_1)^{-1}p\Bigr)
\nonumber\\
&= &\lambda_p(h_2+g_2h_1)S(g_2g_1)\psi\Bigl((g_2g_1)^{-1}p\Bigr)
\nonumber\\
&= &\Bigl([(h_2,g_2)(h_1,g_1)]\psi\Bigr)(p).\nonumber
\end{eqnarray} 
Restricting to $G$, viewed as a subgroup of $H\rtimes_{\alpha}G$, gives a representation of $G$ on $V$ specified by
\begin{equation}\label{E:repGSV}
\bigl({\tilde S}(g)\psi\bigr)(p)=(g\cdot\psi)(p)=S(g)\psi(g^{-1}p).
\end{equation}
We fix a basepoint
$$p_0\in M.$$
We now construct a category ${\mathbf  \Psi}_{p_0}$ whose object set is  $  V$ and whose morphisms (other than the identities) are of the form $(h,g; \psi)\in  (H\rtimes_{\alpha} G)\times \Psi$, with source and target given by 
\begin{equation}\label{E:stargPsiS}
s(h,g, \psi)=  \psi(gp_0) \qquad\hbox{and}\qquad t(h,g, \psi)=  \lambda_{gp_0}(h)\psi(gp_0).
\end{equation}
 We define composition of morphisms by

\begin{equation}\label{E:compos21}
(h_2, g_2,  \psi_2)\circ (h_1, g_1,  \psi_1)=(g_1g_2^{-1}h_2+h_1, g_1,  \psi_1),
\end{equation}
assuming that the source of $(h_2, g_2,  \psi_2)$ is the target of $(h_1, g_1,  \psi_1)$. The target of the morphism on the right in (\ref{E:compos21}) is
\begin{eqnarray*}
\lambda_{g_1p_0}\bigl(g_1g_2^{-1}h_2+h_1\bigr)\psi_1(g_1p_0) &=&\lambda_{g_2p_0}(h_2)\lambda_{g_1p_0}(h_1)\psi_1(g_1p_0)\\
&=&\lambda_{g_2p_0}(h_2)\psi_2(g_2p_0)\\
&=& t(h_2, g_2,  \psi_2),
\end{eqnarray*}
as it should be. 
For associativity we compute:
\begin{eqnarray}\label{E:assoccomp}
(h_3, g_3,  \psi_3)\circ \bigl((h_2, g_2,  \psi_2)\circ (h_1, g_1,  \psi_1)\bigr)&= &(h_3, g_3,  \psi_3)\circ (g_1g_2^{-1}h_2+h_1, g_1,  \psi_1)\\
&= &(g_1g_3^{-1}h_3+g_1g_2^{-1}h_2+h_1,g_1,  \psi_1), \nonumber
\end{eqnarray}
and
\begin{eqnarray*} 
\bigl((h_3, g_3,  \psi_3)\circ  (h_2, g_2,  \psi_2)\bigr)\circ (h_1, g_1,  \psi_1) &=  &(g_2g_3^{-1}h_3+h_2, g_2,  \psi_2)\circ (h_1,g_1,  \psi_1)\\
&=&\bigl(g_1g_2^{-1}(g_2g_3^{-1}h_3+h_2)+h_1, g_1,  \psi_1\bigr),
\end{eqnarray*}
which agrees with (\ref{E:assoccomp}).

We introduce an identity morphism for each object by fiat. This  makes ${\mathbf \Psi}_{p_0}$ a category. We form the semidirect product
\begin{equation}
{\mathbf\Psi}_1=  {\mathbf V}\rtimes_{\eta_1}{\mathbf \Psi}_{p_0},
\end{equation}
where
\begin{equation}\label{E:eta1}
\eta_1: \Mor({\mathbf V})\times \Mor({\mathbf \Psi}_{p_0})\to  \Mor({\mathbf V}):(v, h,g,  \psi)\mapsto  \lambda_{gp_0}(h)v, 
\end{equation}
with $p_0\in M$ being the fixed basepoint we have used above;  we specify the action on identity morphisms by
\begin{equation}\label{E:eta1v1wv}
\eta_1(v,1_w)=v
\end{equation}
for all $v, w\in V$. Before verifying that $\eta_1$ has the required properties for an $\eta$-map let us state the nature of the category ${\mathbf\Psi}_1$. 
 The object and morphism sets for ${\mathbf\Psi}_1$  are
 \begin{eqnarray}\label{E:objmorPsi1}
 \Obj({\mathbf\Psi}_1) &= &\{V\}\times V\simeq V.\\
 \Mor({\mathbf \Psi}_1)' &=& V\times (H\times G\times \Psi),\nonumber
 \end{eqnarray}
 where in the second line we have left out the identity morphisms.
Source and targets are given by
 \begin{eqnarray}\label{E:stPsi1}
 s(v; h,g, \psi) &=& s(h,g, \psi)=\psi(gp_0)\\
 s(v;1_w) &=&w \nonumber\\
  t(v; h,g, \psi) &=& t(h,g, \psi)=\lambda_{gp_0}(h)\psi(gp_0)  \nonumber\\
   t(v;1_w) &=&w. \nonumber
 \end{eqnarray}
Composition of morphisms is given by
 \begin{equation}\label{E:Psi1comp}
 (v_2;f_2)\circ(v_1;f_1)= \bigl(\eta_1(v_2,f_1)+v_1; f_2\circ f_1\bigr).
\end{equation}
More explicitly,
 \begin{eqnarray}\label{E:Psi1comp2}
 (v_2; h_2,g_2, \psi_2)\circ (v_1; h_1,g_1, \psi_1)&=&\left(\lambda_{g_1p_0}(h_1)v_2+v_1; g_1g_2^{-1}h_2+h_1, g_1, \psi_1 \right) \\
  (v_2;f_2)\circ (v; 1_w) &=&(v_2+v; f_2) \nonumber \\
  (v; 1_w) \circ (v_1;f_1) &=&\bigl(\eta(v,f_1)+v_1; f_1\bigr),  \nonumber
 \end{eqnarray}
 assuming that the compositions on the left have appropriate  matching source and targets.

 We turn now to verifying that $\eta_1$ satisfies the required properties. We have
\begin{eqnarray*}
\eta_1\Bigl(\eta_1\bigl(v, (h_2,g_2,  \psi_2)\bigr) , (h_1,g_1,  \psi_1)\Bigr) &=&\eta_1\Bigl(\lambda_{g_2p_0}(h_2)v, (h_1,g_1,  \psi_1)\Bigr)  \\
&=&\lambda_{g_1p_0}(h_1)\Bigr(\bigl(\lambda_{g_2p_0}(h_2)v\bigl)\Bigr)  \\
&=&\lambda_{g_1p_0}(h_1+g_1g_2^{-1}h_2)v
\end{eqnarray*}
which agrees with

 \begin{eqnarray*}
\eta_1\bigl(v, (h_2,g_2,  \psi_2)\circ (h_1,g_1,  \psi_1)\bigr) &=& \eta_1\bigl(v, (g_1g_2^{-1}h_2+h_1, g_1, \psi_1)\bigr)\\
&=&\lambda_{g_1p_0}(h_1+g_1g_2^{-1}h_2)v.
\end{eqnarray*} 
Next,
 \begin{eqnarray*}
\eta_1\bigl(v_2, (h,g,  \psi)\bigr) \circ\eta_1\bigl(v_1, (h,g,  \psi)\bigr) &=&\lambda_{gp_0}(h)v_2+\lambda_{gp_0}(h)v_1\\
&=&\lambda_{gp_0}(h)(v_2+v_1)\\
&=&\eta_1\left(v_2\circ v_1, (h,g,  \psi)\right).
\end{eqnarray*} 
As for  the behavior of $\eta_1$ with respect to identity morphisms we have, in addition to (\ref{E:eta1v1wv}),
$$\eta_1(0, f)= 0$$
for all $f\in\Mor({\mathbf \Psi}_{p_0})$.

  We now work with the categorical group whose object set is $G$ and whose morphism set is $V\rtimes_S G$, with source, target, and $\alpha$-map given by
  \begin{eqnarray*} 
  s(v,g)=t(v,g) &=&g\\
  \alpha(g)v&=&S(g)v.
  \end{eqnarray*} 
  Let us denote this categorical group by
 \begin{equation}
 \mbV\times_{S} \mbg_d,
 \end{equation}
 where $\mbg_d$ is the discrete category with  object set  $G$. Product and composition are given by
 \begin{eqnarray}\label{E:compVeta0G}    (v_2;  g)\circ (v_1;  g) & =& (v_2+v_1;  ,g)\\
 (v_2,g_2)(v_1,g_1)  &=& \bigl(v_2+g_2v_1;  g_2g_1\bigr).\nonumber
   \end{eqnarray}
(Note that we are not viewing $ \mbV\times_{S} \mbg_d$ as an $\eta$-twisted product.)
 We have
 $$[(v'_2,g_2)(v'_1,g_1)]\circ [(v_2,g_2)(v_1,g_1)]=(v'_2+g_2v'_1+v_2+g_2v_1,g_2g_1) $$
 and
 $$(v'_2+v_2,g_2)(v'_1+v_1,g_1)=(v'_2+v_2+g_2v'_1+g_2v_1, g_2g_1)$$

 We define
 \begin{eqnarray}\label{E:eta2Psi}
 \eta_2:  \Mor(\mbV\times_{S} \mbg_d)\times \Mor\left({\mathbf\Psi}_1\right) &\to  &\Mor(\mbV\times_{S} \mbg_d)\\
   \left((v',g'), (v;h,g, \psi)\right) &\mapsto & \bigl(\lambda_{gp_0}(h)v', g'\bigr)\nonumber
 \end{eqnarray}
 and
 \begin{equation}\label{E:eta2id}
 \eta_2\bigl((v',g'), 1_a\bigr)=(v',g').
 \end{equation}
 Before proceeding further let us note the behavior of $\eta_2$ with respect to the other identity morphisms $(0,g')$:
 \begin{equation}
 \eta_2\bigl((0,g'), f\bigr)=(0,g'),
 \end{equation}
 for all $g'\in G$ and all $f\in \Mor({\mathbf\Psi}_1)$.
 
 We now check that $\eta_2$ satisfies the conditions  we have required of $\eta$-maps.  Let
 $$f_1=(v_1;h_1,g_1, \psi_1)\quad\hbox{and}\qquad f_2=(v_2; h_2,g_2, \psi_2). $$
 Then
 \begin{equation}
 f_2\circ f_1 = \bigl(\lambda_{g_1p_0}(h_1)v_2+v_1; g_1g_2^{-1}h_2+h_1, g_1, \psi_1).
 \end{equation}
We have
 \begin{eqnarray*}
 \eta_2\bigl(\eta_2\left((v,g),f_2\right), f_1\bigr) 
 & =&\eta_2\left(\bigl(\lambda_{g_2p_0}(h_2)v,g\bigr), (v_1;h_1,g_1, \psi_1)\right)\\
 &=&\bigl(\lambda_{g_1p_0}(h_1)\lambda_{g_2p_0}(h_2)v, g\bigr)=\bigl(\lambda_{g_1p_0}(h_1+g_1g_2^{-1}h_2)v,  g\bigr)\\
 &=&\eta_2\Bigl((v,g), \bigl(\lambda_{g_1p_0}(h_1)v_2+v_1; g_1g_2^{-1}h_2+h_1, g_1, \psi_1\bigr)\Bigr)\\
 &=&\eta_2\bigl((v,g),\, f_2\circ f_1\bigr).
 \end{eqnarray*}

 Next, working with $g'_2=g'_1\in G$, we have
 \begin{eqnarray*}
 \eta_2\bigl((v'_2,g'_2), (v; h,g, \psi )\bigr) \circ \eta_2\bigl(v'_1,g'_1),  (v;h,g, \psi )\bigr)  &=&\bigl(\lambda_{gp_0}(h)v'_2, g'_2\bigr)\circ \bigl(\lambda_{gp_0}(h)v'_1, g'_1\bigr)\\
& =&\bigl(\lambda_{gp_0}(h)(v'_2+v'_1), g'_1\bigr)\\
&=&\eta_2\bigl((v'_2,g'_2)\circ(v'_1,g'_1),  (v; h,g, \psi)\bigr) 
 \end{eqnarray*}
 
 The condition (\ref{E:etak1g}) also holds for $\eta_2$.
 
We are going to construct a representation $\rho$ of the categorical group $\mbV\times_S\mbg_d$ on ${\mathbf\Psi}_1$ with the twist $\eta_2$:
$$\rho: (\mbV\times_S\mbg_d)\rtimes _{\eta_2}{\mathbf\Psi}_1\to {\mathbf\Psi}_1.$$
Recall that
\begin{eqnarray*}
\Obj(\mbV\times_S\mbg_d) &=&G\\
\Obj({\mathbf\Psi}_1) &=&V.
\end{eqnarray*}

At the level of objects we specify $\rho$ by
\begin{equation}\label{E:rhoobj}
\rho: G\times V\to V:\left(g,v\right)\mapsto gv.
\end{equation}
The morphisms of the two categories are
\begin{eqnarray*} \Mor(\mbV\times_S\mbg_d) &=&V\rtimes_S G\\
\Mor({\mathbf\Psi}_1)' &=&V\times (H\times G\times \Psi),
\end{eqnarray*}
where, as before, we have omitted the identity morphisms in the second line.
At the level of morphisms we define first
 \begin{eqnarray}\label{E:rhomor}  \rho: (V\rtimes_{S}G)\times \Mor({\mathbf \Psi}_1) & \to  &\Mor({\mathbf \Psi}_1) \\
  \left((v',g'), (v;h,g, \psi)\right) &\mapsto& \bigl( g'v-v'; \,   g'h,  \, g'g, g'\cdot\psi\bigr),  \nonumber
 \end{eqnarray}

 where
 \begin{equation}\label{E:gdotpsi}
 g'\cdot\psi ={\tilde S}(g)\psi,
 \end{equation}
 as given in (\ref{E:repGSV}). Next for the action on   morphisms  of the type $(v,1_w)$, where $v,w\in V$, we set
 \begin{equation}\label{E:rho1v}
 \rho(v',g')(v,1_w)=(g'v-v',1_{g'w}) 
 \end{equation}
 for all $(v',g')\in V\times G$.
 We will use the notation
 $$ax\stackrel{\rm def}=\rho(a, x),$$
 for $a$ an object (or morphism) of ${\mathbf V}\times_S{\mbg}_d$ and $x$ an object (or morphism) of ${\mathbf \Psi}_1$.
 
 For $\Mor(\mbV\rtimes_S\mbg_d)$ source and target are the same for each morphism $(v, g)$, both being $g$. The source  of $\bigl( g'v-v';    g'h,  g'g, g'\cdot\psi\bigr)$ is 
 \begin{eqnarray*}
 s\bigl( g'v-v';  g'h, \, g'g,\, g'\cdot\psi\bigr)&=&(g'\cdot\psi)(g'gp_0)= g'\psi({g'}^{-1}g'gp_0) \\ 
  &=&s(v',g')s(v; h,g, \psi).
 \end{eqnarray*}

 For targets we have
   \begin{eqnarray*}
 t\bigl( g'v-v';   g'h,  g'g, g'\cdot\psi\bigr)&=& g' \lambda_{gp_0}(h)\psi (gp_0) \\
 &=&\lambda_{g'gp_0}(g'h)(g'\cdot\psi)(g'gp_0)\\
 &=&\lambda_{gp_0}(h)g'\psi(gp_0)\\
 &=&t(v',g')t(v;h,g, \psi).
 \end{eqnarray*}

 Applying $(v'',g'')$ after applying $(v',g')$  to (\ref{E:rhomor}) we obtain
\begin{equation}\label{E:vprimgprim}
(v'',g'')\bigl[(v',g')  (v;h,g, \psi)\bigr]=  (g''g'v-g''v'-v'';  g''g'h, g''g'g, \, g''g'\cdot\psi),
\end{equation}
 which is the same as applying
 $$(v'',g'')(v',g')=(v''+g''v',g''g')$$
 to $(v; h,g, \psi)$. The same holds for $\bigl[(v'',g'')(v',g')\bigr](v,1_w)$:
 \begin{eqnarray*}
\bigl[(v'',g'')(v',g')\bigr](v,1_w) &=& (v''+g''v',g''g')(v,1_w) \\
 &=&\Bigl(g''g'v-v''-g''v', 1_{g''g'w}\Bigr)\\
 &&\\
  (v'',g'')[(v',g') (v,1_w)] &=& (v'',g'')(g'v-v',1_{g'w}) \\
&=&\Bigl(g''g'v-g''v'-v'', 1_{g''g'w}\Bigr)
 \end{eqnarray*}

 We continue with the notation
 $$f_1=(v_1;h_1,g_1, \psi_1)\quad\hbox{and}\qquad f_2=(v_2; h_2,g_2, \psi_2). $$
 
 Then
\begin{eqnarray}
\label{E:rhov2g2f2} \rho(v'_2,g'_2)f_2\circ \rho(v'_1, g'_1)f_1 &  &\\
&&\hskip -1in = (g'_2v_2-v'_2; g'_2h_2, g'_2g_2, \, g'_2\cdot\psi_2) \circ(g'_1v_1-v'_1; g'_1h_1, g'_1g_1, \,g'_1\cdot\psi_1)\nonumber \\
 &\hskip -2in =&\hskip -1in\Bigl(\lambda_{g'_1g_1p_0}(g'_1h_1)(g'_2v_2-v'_2) + g'_1v_1-v'_1; \Bigr.\nonumber\\
& &\Bigl.  (g'_1g_1)(g'_2g_2)^{-1}g'_2h_2+g'_1h_1,  g'_1g_1, \,g'_1\cdot\psi_1\Bigr)\nonumber
\end{eqnarray}

   Let us compare this with 
 \begin{eqnarray*}
\Bigl(\eta_2\left( (v'_2,g'_2), f_1\right)\circ (v'_1, g'_1)\Bigr)\cdot \Bigl((v_2,f_2)\circ (v_1,f_1)\Bigr)&& \\
 &&\hskip -2in =  \Bigl(\lambda_{g_1p_0}(h_1)v'_2+v'_1; g'_1\Bigr)  \cdot\Bigl(\lambda_{g_1p_0}(h_1)v_2+v_1, g_1g_2^{-1}h_2+h_1,\, g_1,\, \psi_1\Bigr)   \\
& &\hskip -2in=\Bigl( g'_1\bigl(\lambda_{g_1p_0}(h_1)v_2+v_1\bigr)-\lambda_{g_1p_0}(h_1)v'_2-v'_1; \Bigr.\\
&&\hskip -.25in\Bigl. g'_1(g_1g_2^{-1}h_2+h_1), g'_1g_1,   \,g'_1\cdot\psi_1\Bigr). 
 \end{eqnarray*}

 This is identical with the last expression in (\ref{E:rhov2g2f2}) except for one term  but that term is also in agreement upon noting that $g'_1=g'_2$ in order for the composition $\eta_2\left( (v'_2,g'_2), f_1\right)\circ (v'_1, g'_1)$ to be defined.
 
 Next we work with (\ref{E:rhov2g2f2}) for the case $f_2=(v;1_w)$ where $v\in V$ and $w\in V=\Obj({\mathbf\Psi}_1)$, and 
 $$f_1=(v_1;F_1)=(v_1; h_1,g_1, \psi_1)$$ 
 as before, with the composite $f_2\circ f_1$ defined.  Then
 \begin{equation}\label{E:f2circf1rhoeta2}
 f_2\circ f_1=\bigl(\lambda_{g_1p_0}(h_1) v+v_1; h_1,g_1, \psi_1\bigr)
 \end{equation}
 then, using the composition law from (\ref{E:Psi1comp}) and the obvious notation $g'_1F_1$, 
 \begin{eqnarray}\label{E:rhov2g2f2a}
\rho(v'_2, g'_2)f_2\circ \rho(v'_1,g'_1)f_1 & = & (g'_2v-v'_2; 1_{g'_2w}) \circ (g'_1v_1-v'_1; g'_1h_1, g'_1g_1, \,g'_1\cdot\psi_1) \\
&=&\bigl(\, \eta_1(g'_2v-v'_2; g'_1h_1,g'_1g_1, \, g'_1\cdot\psi_1)+g'_1v_1-v'_1;g'_1F_1 \bigr) \nonumber\\
&=&\bigl(\lambda_{g'_1g_1p_0}(g'_1h_1)(g'_2v-v'_2) +g'_1v_1-v'_1;g'_1F_1  \bigr),
 \nonumber
\end{eqnarray}
and, on the other hand,
\begin{eqnarray}\label{E:rhov2g2f2b}
\rho\bigl(\eta_2(v'_2, g'_2; f_1)\circ (v'_1,g'_1)\bigr)(f_2\circ f_1) & =&\rho\bigl(\lambda_{g_1p_0}(h_1)v'_2+v'_1,g'_1)(f_2\circ f_1)  \\
&=& \bigl(g'_1\lambda_{g_1p_0}(h_1) v+g'_1v_1-\lambda_{g_1p_0}(h_1)v'_2-v'_1; \, g'_1F_1 \bigr),\nonumber
  \end{eqnarray}
which agrees with the last expression in (\ref{E:rhov2g2f2a}) because $g'_1=g'_2$, since the composition $(v'_2,g'_2)\circ (v'_1,g'_1)$ is needed to be well-defined. We omit the lengthy but entirely routine verifications of other similar conditions for the cases involving the identity morphisms.

  The twist $\eta_2$ satisfies the invariance relation  (\ref{E:etalkf}) with respect to the representation $\rho$, as may be verified by computation.

 \section{Twisted representations on vector categories}\label{s:irred}
 
 In this section we study actions on categories that have some vector space structure built in.   
 
 By a {\em vector category} we  mean a category whose object and morphism sets are equipped with vector space structures  in such a way that the source and target maps are linear,  as is the identity assigning map $x\mapsto 1_x$.     (This is different from the notion of `linear category' discussed in     \cite{BaezBFW}.)  

By a {\em representation}  $\rho$ of $\mbG$ on a vector category $\mbV$, associated with an $\eta$-map $\eta: \Mor(\mbg)\times \Mor(\mbV)\to\Mor(\mbg)$, we mean an action $\rho$, associated with $\eta$, for which the map
$$\rho(g): \Obj(\mbV)\to\Obj(\mbV): v\mapsto \rho(g,v)$$
and the map
$$\rho(k): \Mor(\mbV)\to\Mor(\mbV): f\mapsto \rho(k,f)$$
are   linear for every object $g$ and every morphism $k$ of $\mbg$.  We shall also call such a representation an {\em $\eta$-twisted representation}, and refer to
$$(\eta, \rho, \mbV)$$
as a {\em twisted representation} of $\mbg$.

 In developing the notion of irreducible representations of categorical groups we need to be careful in handling the role of morphisms that have $0$ as both target and source. This issue is absent in traditional representation theory.

   Let us define the {\em trivial subcategory}  ${\mathbf O}_{\mbV}$ of  a vector category $\mbV$ to be the vector category whose only object element is $0$ and the morphism set is  $ \Hom(0, 0) $:
$$\Obj\bigl({\mathbf O}_{\mbV}\bigr)=\{0\}\qquad \Mor\bigl({\mathbf O}_{\mbV}\bigr)=\Hom(0, 0).$$
 We define a {\em proper vector subcategory} of  $\mathbf V$ to be a subcategory other than $\mathbf V$, ${\mathbf O}_{\mbV}$ and ${\mathbf O}$ (the trivial subcategory with only the $0$ object and its identity morphism).

It is readily verified that an $\eta$-twisted representation $\rho$ of   a categorical group $\mbg$ on $\mbV$,  restricts  to yield a representation of $\mbg$ on ${\mathbf O}_{\mbV}$. With this in mind, we formulate the following definition of reducible representations.
Let $\mbg$ be a categorical group and $\mbV$   a vector category. Let $\rho$ be an $\eta$-twisted representation of $\mbg$ on $\mbV$. We say that $\rho$ is {\em reducible} if there exists a proper  vector subcategory $\mbU$ of $\mbV$ such that the restriction of $(\eta, \rho)$ to $\mbU$ defines 
a representation of $\mbg$ on $\mbU$.  In this case,
\Beq\label{def:red}
\rho|_{\mbU}:\mbg \ltimes_{\eta|_{\mbU}} \mbU\longrightarrow \mbU
\Eeq
is a representation, where $\eta|_{\mbU}$ denotes the restriction of $\eta$ to $\Mor \bigl(\mbU\bigr)\subset \Mor \bigl(\mbV\bigr)$,
$$\eta|_{\mbU}:\Mor \bigl(\mbg\bigr)\times \Mor \bigl(\mbU\bigr)\longrightarrow \Mor \bigl(\mbg\bigr).$$
The representation $\rho$ is  {\em irreducible} if it is not reducible.

We say that a functor  $$\mbF:\mbV\to\mbW,$$
where $\mbV$ and $\mbW$ are vector categories, is {\em linear} if it  induces linear maps both on the object spaces and on the morphism spaces.
A functor $\mbF:\mbV\longrightarrow \mbW$ is a {\em zero functor} if ${\rm Im}(\mbF)\subset{\mathbf O}_{\mbW}$.  Note that the
image of a zero functor at the level of morphisms is not necessarily $\{0\}$.

Suppose $\mbV$ is a vector category, and $\mbV_0$ is a subcategory which is a subspace of $\mbV$ both for objects and for morphisms. Then we would like to define a vector category
$$\mbV/\mbV_0$$
whose object and morphism spaces are the corresponding quotient spaces. However, composition of morphisms would fail to be well-defined since there is no reason why $(f+f_0)\circ(h+h_0)$ should equal  $f\circ h$ modulo $\Mor(\mbV_0)$ when $f_0,h_0\in\Mor(\mbV_0)$. However, as we shall see in our next result the quotient category is well-defined and useful in some cases.

We turn now to Schur's lemma for   irreducible representations of categorical groups on vector categories. There are subtle but important differences between this result and the traditional Schur's lemma for group representations, where the interaction between morphisms and objects is absent.

\begin{theorem}\label{T:schur}
Suppose $(\eta_1, \rho_1, \mbV_1)$ and  $(\eta_2, \rho_2, \mbV_2)$ are irreducible twisted representations of a categorical group 
$\mbg$ on   vector categories $\mbV_1$ and $\mbV_2$, respectively. 
 Let    $\mbF:(\eta_1, \rho_1, \mbV_1)  \longrightarrow (\eta_2, \rho_2, \mbV_2)$  be  
a morphism of representations, given by a linear functor $\mbF:\mbV_1\longrightarrow \mbV_2$.  Then one of the following hold: 
\begin{enumerate}
\item $\mbF$  is a zero functor,\ 
{\hskip 1 cm}or
\item  there is a well-defined quotient category $\mbV_1/\mbO_{\mbV_1}$ and  $\mbF$ induces a functor  
$$\mbV_1/\mbO_{\mbV_1}\longrightarrow  \mbV_2$$
 that is an isomorphism both on objects and on morphisms, or
\item  $\mbF$  is an isomorphism.
\end{enumerate}
\end{theorem}

\pf
The functor $\mbF$ induces linear maps 
$$ \Obj(\mbV_1)\to\Obj(\mbV_2)\qquad\hbox{and}\qquad \Mor(\mbV_1)\to\Mor(\mbV_2).$$

Consider   the  category
\begin{equation}\label{E:defmbL}
\mbL= \ker\mbF ,
\end{equation}
whose objects are those objects in $\mbV_1$ that are mapped by $  \mbF$ to $0\in\Obj(\mbV_2)$ and whose morphisms are those morphisms in $\mbV_1$ that are mapped to the morphism $0\in\Mor(\mbV_2)$.   Let us verify that $\mbL$ is in fact a category: composites of morphisms in $\mbL$ are morphisms in $\mbL$ and for each object $a\in\Obj(\mbL)$ the identity morphism $1_a$ is in $\Mor(\mbL)$.  Since, by definition of a vector category, the mapping $x\mapsto 1_x$ is linear, we have $0=1_0$, the identity map at the zero object $0\in\Obj(\mbV_2)$; so if $f:a\to b$ and $g:b\to c$ are morphisms in  $\mbL$ then
\begin{equation}\label{E:mbLcomp}
{\mathbf F}(g\circ f)={\mathbf F}(g)\circ {\mathbf F}(f)=1_0\circ 1_0=1_0=0,
\end{equation}
and so $g\circ f\in\Mor(\mbL)$. Next, if $a\in\Obj(\mbL)$ then 
\begin{equation}\label{E:mbL1a}
{\mathbf F}(1_a)=1_{{\mathbf F}(a)}=1_0=0,
\end{equation}
and so $1_a\in\Mor(\mbL)$.

The map $\eta_1$ restricts to a map
$$\mbg\times\mbL\to\mbg.$$
We will now show that $\rho$ restricts to a functor
$$\mbg\rtimes_{{\eta}_1}\mbL\to\mbL.$$
Since $\mbF$ is a morphism of representations we have
\begin{equation}\label{E:Frho12L}
\mbF\circ\rho_1=\rho_2\circ ({\rm Id}_{\mbg}\times\mbF)
\end{equation}
(by the commuting diagram (\ref{D:semimap})). Applying this to the object $(g,v)\in \Obj(\mbG)\times\Obj(\mbV_1)$, where $\mbF(v)=0$, we see that
$$\mbF\Bigl(\rho_1(g,v)\Bigr)=\rho_2(g,0)=0,$$
and so the object space of $\mbL$ is stable under the action of $\Obj(\mbg)$ through  the representation $\rho_1$. Similarly,  $\rho_1(k)$ maps $\Mor(\mbL)$ into itself, for every $k\in\Mor(\mbG)$. Thus, $\rho_1$ {\em restricts to an action of $\mbG$  on the subcategory} $\mbL$. Since $\rho_1$ is irreducible it follows that $\mbL $ is  $\mbV_1$ or   $\mbO_{\mbV_1}$ or the zero category $\mbO$. In the first case the functor $\mbF$ would be a zero functor. 

Let us assume then that $\mbF$ is not a zero functor. In this case $\mbL$ is the trivial subcategory $\mbO$ in $\mbV_1$ or $\mbL=\mbO_{\mbV_1}$:
\begin{equation}\label{E:L0V1}
\mbL=\mbO \quad\hbox{or}\quad \mbL=\mbO_{\mbV_1}.
\end{equation}

Let $f, h\in \Mor(\mbV_1)$ for which the composite $f\circ h$ is defined. Let $f_0,h_0\in \Mor(\mbL)$. Then $f+f_0$ has the same source and target as $f$, and the same holds for $h+h_0$ and $h$, and so $(f+f_0)\circ (h+h_0)$ is defined. Moreover,
\begin{eqnarray*}
\mbF\bigl((f+f_0)\circ (h+h_0)\bigr) &=&\mbF(f+f_0)\circ \mbF(h+h_0)\\
&=&\mbF(f)\circ\mbF(h) \qquad\hbox{(because $f_0, h_0\in
\Mor(\mbL)$)} \\
&= &\mbF(f\circ h),
\end{eqnarray*}
and so
$$(f+f_0)\circ (h+h_0) -f\circ h\in \Mor(\mbL).$$
Hence the composition of morphisms in $\mbV_1$ descends to a well-defined composition law on
$$\Mor(\mbV_1)/\Mor(\mbL).$$
Furthermore, as is seen more easily, the source and target maps on $\mbV_1$ induce  well-defined corresponding maps  
$$\Mor(\mbV_1)/\Mor(\mbL)\to \Obj(\mbV_1)=\Obj(\mbV_1)/\Obj(\mbL).$$
Thus we obtain a well-defined vector category
$$\mbV_1/\mbL $$
and a functor
\begin{equation}
\mbF_*:\mbV_1/\mbL\to \mbV_2:\begin{cases} v\mapsto \mbF(v) &\hbox{for $v\in \Obj(\mbV_1)$}\\
f+\Mor(\mbL)\mapsto \mbF(f)&\hbox{for $f\in \Mor(\mbV_1)$}
\end{cases},
\end{equation}
whose image is the same as ${\rm Im}(\mbF)$.  The latter is a $\rho_2$-invariant subcategory of $\mbV_2$ and therefore is either $\mbO_{\mbV_2}$ or $\mbV_2$ itself.
Since  we have assumed that the original functor $\mbF$  is not  a zero functor, its image is, in fact, all of $\mbV_2$.  To conclude let us note that $\mbV_1/\mbL $ is just $\mbV_1$ if $\mbL={\mathbf O}$, and it is $\mbV_1/\mbO_{\mbV_1} $ if $\mbL=\mbO_{\mbV_1}$.
\epf

\section{Representations on categorical vector spaces}\label{S:exm}

The narrowest type of category of interest on which one might consider representations of categorical groups is a categorical vector space. By this we mean a category $\mbV$ for which both object set and morphism set are equipped with vector space structures over a given field $\mbf$, the target and source maps are linear as is the identity assigning map $x\mapsto 1_x$, and the operations of addition  and multiplication by scalar are functorial in the sense that if $f, f', g,g'\in \Mor(\mbV)$, for which $f'\circ f$ and $g'\circ g$ are defined then 
\begin{eqnarray}\label{E:catvec}
(f'+g')\circ (f+g) &= &f'\circ f+g'\circ g\\
(\lambda f')\circ  f &=& \lambda (f'\circ f) =f'\circ (\lambda f)\nonumber
 \end{eqnarray}
for all $\lambda\in{\mathbb F}$.

Thus a categorical vector space is a very special type of vector category:  composition and the vector operations are intertwined through the relation (\ref{E:catvec}).  In this section we work out the structure of categorical vector spaces and the nature of (untwisted) representations on them.

Our first result, formulated in a way that focuses only on the additive abelian group structure of the object and morphism groups of a categorical vector space, demonstrates how special the structure of a categorical vector space is.

\begin{prop}\label{P:abelian2grp} Suppose $\mbV$ is a categorical group for which both object and morphism groups are abelian. Then $\Mor(\mbV)$ is isomorphic to the ordinary product:
\begin{equation}\label{E:MorVW}\Mor(\mbV)\simeq  W\times V,
\end{equation}
where $V$ is the vector space $\Obj(\mbV)$ and $W=\ker s$ is the kernel of the source homomorphism $s:\Mor(\mbV)\to\Obj(\mbV)$. The source and target maps are given on $W\times V$ by
\begin{equation}\label{E:swvtwv}
s(w,v)=v\qquad\hbox{and}\qquad t(w,v)=\tau_0(w)+v,
\end{equation}
where 
\begin{equation}\label{E:swvtwvtau0}
\tau_0:W\to V
\end{equation}
 is   a homomorphism. Moreover,   composition of morphisms corresponds on $W\times V$ to the operation
\begin{equation}\label{E:w2v2w1v1comp}
(w_2,v_2)\circ (w_1,v_1)=(w_2+w_1, v_1),
\end{equation}
defined when $v_2={\tau_0}(w_1)+v_1$. The identity morphism at any $v\in V$ is given by $(0,v)$, and the inverse of a morphism $(w,v)$ is $\bigl(-w,\tau_0(w)+v\bigr)$.
\end{prop}

\pf  As discussed earlier in the context of (\ref{E:morcatcross}), 
\begin{equation}\label{E:MrVVW}
\Mor(\mbV)\simeq W\rtimes_{\alpha}V
\end{equation}
 where
$$\alpha:V\times W\to W: (v,w)\to \alpha(v)(w)=\alpha_v(w)$$
is a mapping for which  each $\alpha_v$ is an automorphism of the group $W$.  Source and target maps are given by
\begin{equation}\label{E:catst}
s(w,v)=v\quad\hbox{and}\quad t(w,v)=\tau_0(w)+v,
\end{equation}
where, as we saw in (\ref{E:thgtau}), $\tau_0$ is the restriction of $t$ to $W$:
\begin{equation}\label{E:tau0def}
\tau_0=t|W.
\end{equation}

The binary operation for the semidirect product $W\rtimes_{\alpha}V$ is given by
\begin{equation}\label{E:w2v2plus}
(w_2,v_2)\oplus (w_1,v_1)=(w_2+\alpha_{v_2}(w_1), v_2+v_1).\end{equation}
Thus $W\rtimes_{\alpha}V$ is abelian if and only if 
$$(w_2+\alpha_{v_2}(w_1), v_2+v_1)=(w_1+\alpha_{v_1}(w_2), v_1+v_2)$$
for all $v_1, v_2\in V$ and $w_1, w_2\in W$. Comparing the first components, and setting $w_2=0$ and $v_1=0$ we have as necessary condition:
$$\alpha_{v}(w)=w $$
for all $v\in V$ and $w\in W$. This means that $\alpha_v$ is the identity map on $W$ for all $v$, which is equivalent to the group operation (\ref{E:w2v2plus}) being given by
\begin{equation}\label{E:w2v2plusabel}
(w_2,v_2)\oplus (w_1,v_1)=(w_2+w_1, v_2+v_1).\end{equation} 
This condition is obviously also sufficient to ensure that $W\rtimes_{\alpha}V$ is abelian. In this case the semidirect product $W\rtimes_{\alpha}V$ is just the ordinary product $W\times V$:
$$W\rtimes_{\alpha}V=W\times V.$$  

We note that $(0,v_2)\circ (w_1,v_1)=(w_1,v_1)$ when $v_2=t(w_1,v_1)$, and
$(w_2, v_2)\circ (0,v_1)=(w_2, v_1)$ when $v_2=t(0,v_1)=v_1$.  It follows then that the identity morphism at any $v$ is $(0,v)$. Next for a morphism $(w_1,v_1)$ we have
$$(-w_1,v'_1)\circ (w_1,v_1)=(0,v_1),$$
where $v'_1=\tau_0(w_1)+v_1$, and 
$$(w_1,v_1)\circ (-w_1, v'_1)=(0,v'_1)$$
because $\tau_0(-w_1)+v'_1=v_1$; thus the inverse of $(w_1,v_1)$ is $(-w_1,v_1)$. \epf

Now we specialize to a categorical vector space $\mbV$. Writing the isomorphism (\ref{E:MrVVW}) as given in (\ref{E:morcatcross})  explicitly, using the additive notation, we have 
\begin{equation}\label{E:MorVWiso}
\Mor(\mbV)\to W\rtimes_{\alpha}V: f\mapsto (f+1_{-s(f)}, s(f)\bigr).
\end{equation}
Linearity of the source map $s$ and of the identity-assigning map $a\mapsto 1_a$ imply then that the isomorphism in (\ref{E:MorVWiso}) is an  {\em  isomorphism of vector spaces.} The mapping ${\tau_0}$  given  in Proposition \ref{P:abelian2grp} above is then a {\em linear } map
$${\tau_0}:W\to V.$$

Consider a representation (without any twist) $\rho$  of $\mbg$ on $\mbV$; we use notation $G$, $H$, and $V$, $W$ as above.  Thus $\mbg$ is associated to the crossed module $(G,H,\alpha,\tau)$. To keep the notation simple, we denote the representation of $K= \Mor(\mbg)$ on $W$ as well as the representation of $G=\Obj(\mbg)$ on $V$ by $\rho$.  We also denote all target and source maps by $t$ and $s$, respectively,   both as maps $\Mor(\mbg)\to\Obj(\mbg)$ and in their more concrete renditions as maps $G\rtimes_{\alpha}H\to G$.
We will write $\rho(k)f$ as $kf$  when the intended meaning is clear.

 Let us first see how $\rho$ introduces a representation of  the group $K$ on the vector space $W$. To this end we observe that when an element $k\in K$ acts on an element of the form $(w,0)\in W\times V$, the source of the resulting element is
 $$s[\rho(k)(w,0)]=\rho\bigl(s(k)\bigr)s(w,0)=\rho\bigl(s(k)\bigr)0=0.$$
This means that the subspace  $W\times\{0\}$ of $W\times V$ is invariant under the representation $\rho$ of $K$. Thus $\rho$ produces a representation, which we   denote by $\rho_0$, of $K$ on $W$:   
\begin{equation}\label{E:defrho}
\rho(k)w= \rho_0(k)w \qquad\hbox{for all $w\in W$ and $k\in K$.}
\end{equation}
  To see how this interacts with the representation $\rho$  of $G=\Obj(\mbg)$ on $V$ we compute 
\begin{equation}
t\bigl(\rho_0(k)w\bigr)=t[\rho(k)w]=\rho\bigl(t(k)\bigr)t(w)=\rho\bigl(t(k)\bigr){\tau_0}(w).
\end{equation}

For the interaction of $\rho_0$ and the composition of morphisms in $K$, let us recall that
$$\rho:\mbg\times\mbV\to\mbV$$
is a functor. Using this and the definition of $\rho_0$ in  (\ref{E:defrho}) we have:
\begin{eqnarray}\label{E:defrhotcomp} 
\bigl(\rho_0(k_2\circ k_1)w,0\bigr) &=&\rho(k_2\circ k_1)(w,0) \\
&=& [\rho(k_2\circ k_1)] [(w,0)\circ (0,0)] \nonumber\\
&=&[\rho(k_2)(w,0) ]\circ[ \rho(k_1)(0,0)]\nonumber\\
&=&\bigl(\rho_0(k_2)  w,0\bigr)\circ (0,0)\nonumber\\
&=&\bigl(\rho_0(k_2)  w,0\bigr)\nonumber
\end{eqnarray}
 for any $k_1, k_2\in K=\Mor(\mbg)$ for which the composite $k_2\circ k_1$ exists and all $w\in W$. This is a very strong restriction on   $\rho_0$. In fact, taking $k_1=k:a\to b$ and $k_2=1_b$, the identity morphism at $b$, we have
\begin{equation}\label{E:rho0k1b}
\rho_0(k)=\rho_0(1_{t(k)}).
\end{equation} 
 Thus $\rho_0$ is completely determined by the representation of $G$ on $W$ given by
\begin{equation}\label{E:rhoGW}
G\to \End(W): b\mapsto \rho_0(1_b).
\end{equation}
 We denote this representation of $G$ by $\rho_0$:
 \begin{equation}\label{E:rho0bdef}
 \rho_0(b)\stackrel{\rm def}{=}\rho_0(1_{ b}).
 \end{equation} 
 
To summarize, we have proved:

\begin{prop}\label{P:repcatvec}
Let $\mbg$ be a categorical group,  with source map $s$ and target map $t$, and let $(G, H,\alpha,\tau)$ be the associated  crossed module.   Suppose $\rho$ is a representation of $\mbg$  on a categorical vector space $\mbV$. Let     $V=\Obj(\mbV)$, $W=\ker s$, and ${\tau_0}:W\to V$ be the restriction of the target map to $W$.  Then there is a representation $\rho_0$ of $G$ on $W$ such that
\begin{equation}
\rho(k)w= \rho_0\bigl(t(k)\bigr)w \qquad\hbox{for all $w\in W$ and $k\in K$,}
\end{equation}
and  $\rho_0$ intertwines with $\rho|G$ through the map ${\tau_0}:W\to V$:
\begin{equation}\label{E:tau0intert}
{\tau_0}\bigl(\rho_0(k)w\bigr)=\rho(t(k)\bigr){\tau_0}(w)\quad\hbox{for all $k\in K$ and $w\in W$.}
\end{equation}

 \end{prop}
 
 The very special conditions that $\rho$ must of necessity satisfy can be  explored further.  Let us recall from (\ref{E:100}) that the zero element in $\Mor(\mbV)$ is $1_0$:
 $$0=1_0,$$
 and that since $\rho$ is a functor it carries identity morphisms to identity morphisms:
 \begin{equation}\label{E:rho1ga}
 \rho(1_g,1_v)=1_{gv}.
 \end{equation}
for all $g\in \Obj(\mbg)$ and $v\in \Obj(\mbV)$. For any morphism $f:v\to w$ in $\Mor(\mbV)$,  $\rho(1_g,f)$ has source $gv$ and target $gw$ while   $\rho(1_g,f^{-1})$ has source $gw$ and target $gv$; then:
\begin{equation}
1_{gv}=\rho(1_g, f^{-1}\circ f) =\rho(1_g\circ 1_g, f^{-1}\circ f)=\rho(1_g,f^{-1})\circ \rho(1_g,f)
\end{equation}
and
   \begin{equation}
1_{gw}=\rho(1_g, f\circ f^{-1}) =\rho(1_g\circ 1_g, f\circ f^{-1})=\rho(1_g,f)\circ \rho(1_g,f^{-1}).
\end{equation}
Hence 
\begin{equation}\label{E:1gfinv}
 1_g f^{-1} =\bigl(1_gf\bigr)^{-1}.
\end{equation}
We can combine this with the very special relation (\ref{E:rho0k1b}
)
$$\rho_0(k)=\rho_0(1_{t(k)})$$
to extract more information on the nature of $\rho$.

 If $f=(w,v)\in W\oplus V\simeq \Mor(\mbV)$ has {\em source $s(f)=v$ in the range of}  $\tau_0$ then we can express $f$ as a   composite using  two elements in  $W$:
 $$f=(f\circ f_0)\circ f_0^{-1},$$
 where $f_0=(w_v,0)\in W$ and $w_v$ being such that $\tau_0(w_v)=v$, and $f\circ f_0$ is also in $W$ because its source is $0$.   Then   
 \begin{equation}
 \rho(k,f)=\rho(k, f\circ f_0)\circ \rho(1_b, f_0^{-1})=\rho(k, f\circ f_0)\circ [\rho(1_b, f_0)]^{-1}
 \end{equation}
 where $b=t(k)$. We recall from (\ref{E:rho0k1b})  that 
\begin{eqnarray*}
\rho(k, f\circ f_0) &=&\rho(1_{t(k)}, f\circ f_0)=\rho_0(b)(f\circ f_0).\\
\rho(1_b,f_0) &=&\rho_0(b)f_0.
\end{eqnarray*}
  Our conclusion then is:
 
{\em the values $\rho(k,f)$, when $s(f)\in t(W)$, are uniquely determined by the  representation $\rho_0$ of $G$ on $W$.}

The value of $\rho(k,f)$ for $f:v \to w$  is also uniquely determined by the restriction of the representation $\rho$ of $G$ to $\Hom(0,0)\subset\Mor(\mbV)$ and the value $\rho(k,f')$ for any one $f'\in \Hom(v,w)$. This is seen by writing $f$ as $f'+(f-f')$ and noting that $f-f'\in \Hom(0,0)$.

Thus categorical vector spaces and representations of categorical groups on them have a very special structure.


\section{Concluding summary}

In this paper we have studied actions of categorical groups on categories. The distinction between such actions and traditional group actions lies in the behavior of compositions of morphisms, a structure not present in ordinary groups. We   have introduced a notion of a twisted action in this context, allowing for a richer range of behaviors of composition in relation to the action. For linear actions on categories with vector space structures we have  proved a version of Schur's lemma and studied representations on categorical vector spaces.    To illustrate the ideas we have developed several examples in detail.

 {\bf{ Acknowledgments.} } We thank Manuel B{\"a}renz for comments that led us to  correct an earlier version of Theorem \ref{T:schur}. Sengupta   acknowledges  (i) research support from   NSA grant H98230-13-1-0210, (ii) the American Institute of Mathematics  for participation in a workshop in April 2013, and (iii) the SN Bose National Centre for its hospitality during visits when part of this work was done.  Chatterjee acknowledges support through a fellowship from the Jacques Hadamard Mathematical Foundation. We are very grateful to the anonymous referee for a thorough reading of this paper and for comments that led to corrections and improvements in presentation.
 
\refs

\bibitem[Abbaspour and Wagemann, 2012]{AbbasWag} Hossein Abbaspour and Friedrich Wagemann, {\em On 2-Holonomy}, (2012) \url{http://arxiv.org/abs/1202.2292} 

\bibitem[Awody, 2006]{Awody} Steve Awody, {\em Category Theory}, Clarendon Press, Oxford (2006) 

\bibitem[Baez, 2002]{Baez} John C. Baez, {\em Higher Yang-Mills theory}, (2002) \url{http://arxiv.org/abs/hep-th/0206130}

\bibitem[Baez et al., 2012] {BaezBFW} John C. Baez, Aristide Baratin, Laurent Freidel, and Derek K. Wise,   {\em Infinite-Dimensional Representations of $2$-Groups}, Mem. Amer. Math. Soc. {\bf 219} (2012), no. 1032

\bibitem[Baez and Schreiber, 2007]{BaezSchr} John C. Baez and Urs Schreiber, {\em Higher gauge theory},  Categories in algebra, geometry and mathematical physics, (eds. A. Davydov et al), Contemp. Math., Amer. Math. Soc., {\bf 431}(2007), 7-30

\bibitem[Baez and Wise, 2012]{BaezWise} John C. Baez and Derek K. Wise, {\em Teleparallel Gravity as a Higher Gauge Theory}, (2012) \url{http://arxiv.org/abs/1204.4339}

\bibitem[Baratin and Wise, 2009]{BarWise} Aristide Baratin and Derek K. Wise, {\em 2-Group Representations for Spin Foams}, (2009) \url{http://arxiv.org/abs/0910.1542}

\bibitem[Barrett and Mackaay, 2006]{BarrMack} John W. Barrett and Marco Mackaay,  {\em Categorical Representations of Categorical Groups}. Theory and Applications of Categories, {\bf 16} (2006), no. 20,  529-557

\bibitem[Bartels, 2004]{Bart} Toby Bartels,  {\em Higher Gauge Theory I: 2-Bundles},  (2004) \url{http://arxiv.org/abs/math/0410328}

 \bibitem[Chatterjee et al., 2010]{CLSpath} Saikat Chatterjee, Amitabha Lahiri, and Ambar  N. Sengupta, {\em Parallel Transport over Path Spaces},  Reviews in Math. Phys., {\bf 9} (2010), 1033-1059

\bibitem[Chatterjee et al., 2014]{CLScb}Saikat Chatterjee, Amitabha Lahiri and Ambar N. Sengupta, {\em Path space Connections and Categorical Geometry},  Journal of Geom. Phys.,  {\bf 75} (2014), 129-161 

\bibitem[Crane and Sheppeard, 2003]{CS} Louis Crane and Marni D. Sheppeard, {\em 2-categorical Poincare Representations and State Sum Applications}, (2003) \url{http://arxiv.org/abs/math/0306440}

\bibitem[Elgueta, 2007]{Elgu} Josep Elgueta, {\em Representation theory of 2-groups on Kapranov and Voevodsky 2-vector spaces},
Adv. Math., {\bf 213} (2007), 53-92

\bibitem[Forrester-Barker, 2002]{ForBar} Magnus Forrester-Barker, {\em Group Objects and Internal Categories}, (2002) 	\url{http:/arXiv:math/0212065v1 }

\bibitem[Kamps et al., 1982]{KPT}  Klaus Heiner Kamps, Dieter Pumpl\"un, and Walter Tholen (eds.), Category Theory (Gummersbach,
1981), no. 962 in Lecture Notes in Mathematics, Berlin, 1982, Springer-Verlag

\bibitem[Kelly, 1982]{Kelly} Gregory M. Kelly, {\em Basic Concepts of Enriched Categories}, Cambridge University Press  (1982)

\bibitem[Mac Lane 1971]{Macl} Saunders Mac Lane, {\em Categories for the Working Mathematician}. Springer-Verlag (1971)

\bibitem[Martins and Picken, 2010]{MarPick} Joao F. Martins and Roger Picken, {\em On two-Dimensional Holonomy}, Trans. Amer. Math. Soc. {\bf 362} (2010), no. 11, 5657-5695

\bibitem[Peiffer, 1949]{Peif} Ren\'ee Peiffer, {\em \"Uber Identit{\"a}ten zwischen Relationen}. Math. Ann. {\bf 121} (1949),  67-99

\bibitem[Pfeiffer, 2003]{Pfeif}Hendryk Pfeiffer, {\em  Higher gauge theory and a non-Abelian generalization of 2-form electrodynamics}, Ann. Physics {\bf 308} (2003), no. 2, 447-477

 \bibitem[Porter, 1981]{Port } Timothy Porter,  {\em Internal categories and crossed modules}, in \cite{KPT},  249-255
 
 \bibitem[Schreiber and Waldorf, 2007]{SW1}Urs Schreiber and Konrad Waldorf, {\em Parallel Transport and Functors}, J. Homotopy Relat. Struct. {\bf 4},  (2009), no. 1, 187-244

\bibitem[Schreiber and Waldorf, 2008]{SW2}Urs Schreiber and Konrad Waldorf, {\em Connections on non-abelian Gerbes and their Holonomy},  Theory Appl. Categ. {\bf 28} (2013), 476-540

\bibitem[Steinberg, 1999]{steinberg} Benjamin Steinberg, {\em Semidirect products of categories and applications}, J. Pure Appl. Algebra {\bf 142} (1999), no. 2, 153-182
 
 \bibitem[Viennot, 2012]{Vien} David Viennot, {\em Non-abelian higher gauge theory and categorical bundle}, (2012) \url{http://arxiv.org/abs/1202.2280}
 
 \bibitem[Whitehead, 1949]{Whit} John Henry C. Whitehead, {\em Combinatorial Homotopy II.} Bull. Amer. Math. Soc., {\bf 55} (1949), no. 5, 453-496
 
 \bibitem[Yetter, 2003]{Yet} David N. Yetter, {\em Measurable Categories}, Appl. Categ. Structures {\bf 13} (2005), no. 5-6, 469-500

\endrefs 
\end{document}